%% file: main.tex
\documentclass[sigconf]{acmart}

\newcommand{\ourtitle}{Algorithms for Hiring and Outsourcing in the \mbox{Online Labor Market}}

\usepackage{amsmath}
\usepackage{mathtools}

\newcommand{\figDirectory}{fig}

\copyrightyear{2018} 
\acmYear{2018} 
\setcopyright{acmcopyright}
\acmConference[KDD 2018]{24th ACM SIGKDD International Conference on Knowledge Discovery \& Data Mining}{August 19--23, 2018}{London, United Kingdom}
\acmBooktitle{KDD 2018: 24th ACM SIGKDD International Conference on Knowledge Discovery \& Data Mining, August 19--23, 2018, London, United Kingdom}
\acmPrice{15.00}
\acmDOI{10.1145/3219819.3220056}
\acmISBN{978-1-4503-5552-0/18/08} 

\usepackage[hang,tight]{subfigure}  
\usepackage{booktabs}               

\usepackage{paralist}               

\usepackage{balance}                
\usepackage{xspace}

\def\pprw{8.5in} \def\pprh{11in}

\newtheorem{problem}{Problem}
\theoremstyle{empty}

\newcommand{\spara}[1]{\smallskip\noindent{\bf #1}}

\newsavebox{\mybox}
\newenvironment{framed}{\noindent\begin{lrbox}{\mybox}\begin{minipage}{.98\columnwidth}}{\end{minipage}\end{lrbox}\fbox{\usebox{\mybox}}}

\newcommand{\Ceil}[1]{\left\lceil#1\right\rceil}

\newcommand{\Card}[1]{\left\lvert#1\right\rvert}

\newcommand{\Def}{\overset{\triangle}=}
\newcommand{\St}{;\:}
\newcommand{\E}[1]{\mathbf{E}\!\left[#1\right]}

\newcommand{\PrimalTime}[1]{\ensuremath{P^{#1}}}

\newcommand{\DualTime}[1]{\ensuremath{D^{#1}}}

\newcommand{\skillset}{\ensuremath{{S}}\xspace}
\newcommand{\taskset}{\ensuremath{\mathcal{J}}\xspace}
\newcommand{\workerset}{\ensuremath{\mathcal{W}}\xspace}
\newcommand{\hiredworkers}{\ensuremath{\mathcal{H}}\xspace}
\newcommand{\nonhiredworkers}{\ensuremath{\mathcal{F}}\xspace}

\newcommand{\numskills}{\ensuremath{m}\xspace}
\newcommand{\numtasks}{\ensuremath{T}\xspace}
\newcommand{\numtasksTotal}{\ensuremath{T^{\ast}}\xspace}
\newcommand{\numworkers}{\ensuremath{n}\xspace}

\newcommand{\task}[1]{\ensuremath{J^{#1}}\xspace}
\newcommand{\atask}{\ensuremath{J}\xspace}

\newcommand{\worker}[1]{\ensuremath{W^{#1}}\xspace}
\newcommand{\aworker}{\ensuremath{W}\xspace}

\newcommand{\askill}{\ensuremath{\ell}\xspace}

\newcommand{\hiringcost}[1]{\ensuremath{C_{#1}}\xspace}
\newcommand{\maxhiringcost}[1]{\ensuremath{C^*}\xspace}
\newcommand{\outsourcecost}[1]{\ensuremath{\lambda_{#1}}\xspace}
\newcommand{\maxoutsourcecost}[1]{\ensuremath{\lambda^*}\xspace}

\newcommand{\salary}[1]{\ensuremath{\sigma_{#1}}\xspace}

\newcommand{\hiringcostwithsalary}[1]{\ensuremath{\widehat{C}_{#1}}\xspace}
\newcommand{\maxhiringcostwithsalary}[1]{\ensuremath{\widehat{C}^*}\xspace}

\newcommand{\hireintervals}{\ensuremath{\eta}\xspace}

\newcommand{\solution}{\ensuremath{\mathcal{Q}}\xspace}

\newcommand{\pool}[1]{\ensuremath{P_{#1}}\xspace}

\newcommand{\intervals}{\ensuremath{\mathcal{I}}\xspace}

\newcommand{\SolutionTFO}{\ensuremath{S^{\text{TFO}}}\xspace}
\newcommand{\CostSolutionTFO}{\ensuremath{C^{\text{TFO}}}\xspace}
\newcommand{\SolutionAlt}{\ensuremath{S^{\text{Alt}}}\xspace}

\newcommand{\SolutionLetter}{\ensuremath{S}}

\newcommand{\SolAltApprox}{\ensuremath{\hat{\SolutionLetter}^{\text{Alt}}}\xspace}
\newcommand{\SolTFOApprox}{\ensuremath{\hat{\SolutionLetter}^{\text{TFO}}}\xspace}

\newcommand{\greedy}{\texttt{Greedy\_SC}\xspace}

\newcommand{\altalgo}{\texttt{Alt-TFO}\xspace}
\newcommand{\price}[1]{\ensuremath{\operatorname{price}}(#1)}
\newcommand{\algouseforfree}{\texttt{LumpSum}\xspace}
\newcommand{\naive}{\texttt{LumpSum-Heuristic}\xspace}
\newcommand{\naivetfo}{\texttt{TFO-Heuristic}\xspace}
\newcommand{\algowithsalary}{\texttt{TFO}\xspace}

\newcommand{\algoAlwaysHire}{\texttt{Always-Hire}\xspace}
\newcommand{\algoAlwaysOutsource}{\texttt{Always-Outsource}\xspace}
\newcommand{\metaalgo}{\texttt{TFO-Adaptive}\xspace}

\newcommand{\skirentalproblem}{\textsc{SkiRental}\xspace}
\newcommand{\setcoverproblem}{\textsc{SetCover}\xspace}
\newcommand{\problemuseforfree}{\textsc{LumpSum}\xspace} 
 
\newcommand{\problemsalary}{\textsc{TFO}\xspace}
\newcommand{\ourproblem}{\textsc{TFO}\xspace} 
\newcommand{\problemalternative}{\textsc{Alt-TFO}\xspace}

\newenvironment{enumerate-algo}%
{\begin{list}{\arabic{enumi}.}%
      {\setlength{\leftmargin}{2.5em}%
       \setlength{\itemsep}{-\parsep}%
       \usecounter{enumi}}%
}{\end{list}}

\begin{document}

\setdefaultleftmargin{1em}{}{}{}{.5em}{.5em}

\fancyhead{}

\title{\ourtitle} 
\titlenote{Because of space limitations, some details and
proofs have been omitted, they will appear in the full version of this
  work.}
\titlenote{The research for this work has been supported by
the Google Focused Research Award ``Algorithms for Large-Scale Data
Analysis,''
the EU FET project MULTIPLEX 317532,
the ERC Advanced Grant 788893 ``Algorithmic and Mechanism Design Research for Online MArkets (AMDROMA),''
NSF grants CAREER 1253393 and IIS 1421759, and
La Caixa project LCF/PR/PR16/11110009.}

\author{Aris Anagnostopoulos}
\affiliation{
  \institution{Sapienza University of Rome}
}

\author{Carlos Castillo}
\affiliation{
  \institution{Universitat Pompeu Fabra}
}
\author{Adriano Fazzone}
\affiliation{
   \institution{Sapienza University of Rome}
}
\author{Stefano Leonardi}
\affiliation{
   \institution{Sapienza University of Rome}
}
\author{Evimaria Terzi}
\affiliation{
   \institution{Boston University}
}

\renewcommand{\shortauthors}{Anagnostopoulos et al.}

\begin{abstract}
\input{abstract}
\end{abstract}

\maketitle

\section{Introduction}\label{sec:intro}
\input{intro}

\section{Preliminaries}\label{sec:preliminaries}
\input{preliminaries}

\section{The {\large{\problemuseforfree}} Problem}\label{sec:useforfree}
\input{simple}

\section{The {\large {\problemsalary}} Problem}\label{sec:salarynew}
\input{salarynew}

\section{Experiments}\label{sec:experiments}
\input{experiments}

\section{Related Work}\label{sec:relwork}
\input{related}

\section{Conclusions}\label{sec:conclusions}
\input{conclusions}

\balance
\bibliographystyle{ACM-Reference-Format}
\bibliography{rent-or-buy}






\end{document}

%% file: abstract.tex
Although freelancing work has grown substantially in recent years, in part facilitated by a number of online labor marketplaces, 
traditional forms of ``in-sourcing'' work continue being the dominant form of employment. 
This means that, at least for the time being, freelancing and salaried employment will continue to co-exist.
%
%
In this paper, we provide algorithms for outsourcing and hiring workers in a general setting, where workers form a team and contribute different skills to perform a task.
We call this model \emph{team formation with outsourcing}.
In our model, tasks arrive in an online fashion: neither the number nor the composition of the tasks are known a-priori. 
At any point in time, there is a team of hired workers who receive a fixed salary independently of the work they perform. 
This team is dynamic: new members can be hired and existing members can be fired, at some cost.
Additionally, some parts of the arriving tasks can be outsourced and thus completed by non-team members, at a premium.
%
Our contribution is an efficient online cost-minimizing algorithm for hiring and
firing team members and outsourcing tasks.
We present theoretical bounds obtained using a primal--dual scheme proving that our algorithms have logarithmic competitive approximation ratio. We complement these results with experiments using semi-synthetic datasets based on actual task requirements and worker skills from three large online labor marketplaces.
%

%% file: intro.tex
Self-employment is an increasing trend; for instance, between 10\% and
20\% of workers in 
developed countries are
self-employed~\cite{oecd2016data}.
This phenomenon can be partially attributed to business downsizing and
employee dissatisfaction,
as well as to the existence of online labor markets (e.g., \emph{Guru.com,} \emph{Freelancer.com}).
This trend has enabled freelancers to
work remotely on specialized tasks, 
and prompted researchers
and practitioners to explore the benefits of outsourcing and
\emph{crowdsourcing}~\cite{howe06,jeppesen10,malone10,surowiecki04,kittur11,retelny14}.


Although crowdsourcing adoption was driven, at least in part, by the
assumption that problems can be decomposed into parts that can be
addressed separately by independent workers, 
crowdsourcing results can be
improved by allowing some degree of collaboration among 
them~\cite{majchrzak2013towards,Riedl2016}.
The idea of combining collaboration
with crowdsourcing has led to research on \emph{team
formation}~\cite{anagnostopoulos2010power,
anagnostopoulos12online, an13finding, dorn10composing,
gajewar12multiskill, golshan14profit, kargar11teamexp,
lappas2009finding, li10team, sozio10community, majumder12capacitated},
in which a common thread is the need for complementary skills, and
problem settings differ in aspects such as objectives (e.g., load
balancing and/or compatibility), constraints (e.g., worker capacity),
and algorithmic set up (online or offline).


\spara{Overview of problem setting and assumptions.}
%
%
We consider tasks that arrive in an online fashion and must be completed
by assigning them to one or more workers, who jointly cover the skills required for each task.
At any point in time, there is a team of hired workers who are paid a salary, independently of the work they perform. 
This team is dynamic: new members can be hired and existing members can be fired.
Hiring and firing workers is expensive, which is why companies routinely keep on the payroll skilled workers even if they are temporarily idle; however, they also seek to maintain ``benching'' to a minimum~\cite{sushma2017waning}.
Outsourcing provides additional flexibility as some parts of the incoming tasks can be completed by non-team members who are outsourced.
In practice, outsourcing involves additional costs such as searching, contracting, communicating with, and managing an expert or specialist external to a company~\cite{barthelemy2001hidden}.

Deciding when to hire, fire, and outsource workers is a difficult online problem with parameters that depend on job market conditions and employment regulations.
Intuitively:
\begin{inparaenum}
 \item if the cost of hiring or firing workers is too high, outsourcing becomes preferable to hiring;
 \item if the cost of outsourcing work relative to salaries of hired workers is too high, hiring becomes preferable to outsourcing; and
 \item if the workload consists of many repetitions of similar tasks, hiring becomes preferable to outsourcing.
\end{inparaenum}

In this paper, we formulate this as an online cost minimization problem,
which we call Team Formation with Outsourcing ({\ourproblem}). 
We formally define this problem in Section~\ref{sec:preliminaries} and solve it in
Sections~\ref{sec:useforfree} and~\ref{sec:salarynew}.
%
%
Despite this being a model and hence not capturing every aspect of employment decisions in a company, we show how it brings formalism to the intuitions we have outlined,
helps understand under which circumstances a combination of hiring and outsourcing can be cost effective,
and motivates experimentation on semi-synthetic data allowing us to cover a broad range of cases, as we show in Section~\ref{sec:experiments}.


\spara{Algorithmic techniques.}
To the best of our knowledge, we are the first to consider this problem and
study some of its variants.
Our problem turns out to be an original generalization  of {\em online set
cover} and  {\em online ski rental}, two of the most paradigmatic
online problems.
In fact \ourproblem has elements that make it more complex;
to solve it, an algorithm has to address its various characteristics:
\begin{inparaenum}[(1)]
 \item it is also online, so decisions should be taken with limited information on the input, but at each step, 
 an entirely new instance of the set-cover problem needs to be solved by using hired and outsourced workers;
 \item hired and outsourced workers collaborate with each other, and this needs to be taken into account; and 
 \item workers can be hired, fired, hired again, and so on, so one has to keep track of their status at every point in time. 
\end{inparaenum}

Several natural approaches inspired by  online algorithms for the problems we mentioned previously,
fail to provide solutions with theoretical guarantees.
Therefore, we consider an approach introduced in the last  years 
for studying complex online problems, the \emph{online primal--dual} scheme~\cite{buchbinder2009design}.  
The idea is to create a sequence of integer programs to model the
online problem by incrementally introducing variables and constraints. We then
consider their linear relaxations and their duals to design an online algorithm and we  analyze 
it by comparing the costs of the primal and the dual programs as they
evolve over time with the arrival of new tasks.  This is a powerful
approach, which has so far been
applied with success to several classical online problems: packing and  covering problems, ski-rental, weighted caching, k-server among others \cite{Bansal13}.
We refer to \cite{buchbinder2009design, Bansal13} for a survey of the  applications of the online primal-dual method.

%
%

Our analysis results in polynomial-time algorithms that have logarithmic competitive approximation ratios.  
This means that despite the fact that our algorithms work in an online
fashion and they do not have any knowledge of the number and the
composition of future tasks, we can guarantee that the cost they will
incur will be, \emph{at every time instance,} only a logarithmic factor worse
than the cost incurred by an optimal algorithm that knows the set of requests a priori.

\spara{Contributions.}
The key contributions of our work are:
\begin{compactitem}


\item We formalize {\ourproblem}: the problem  
of designing an online cost-minimizing algorithm for hiring, firing and
outsourcing. 

\item We design efficient and effective approximation algorithms for  {\ourproblem} using an online primal--dual scheme, and provide approximation guarantees on their performance.


\item We experiment on semi-synthetic data based on actual task requirements and worker skills from three large online labor marketplaces, testing algorithms under a broad range of conditions.

\item We provide experimental evidence of the quality of the performance of online primal--dual algorithms for a complex real-world problem. 
Prior work has performed  theoretical analysis mostly for classical or practically motivated  online problems \cite{Buchbinder:2007, Devanur:2009}.  To the best of our knowledge, the empirical validation was previously addressed only for the Adwords matching problem \cite{ho2012online}. We demonstrate that such
approaches, even though they are based on heavy theoretical machinery,
  can be easily implemented and are efficient in practice.
\end{compactitem}

%% file: preliminaries.tex
In this section, we formally describe our setting and problem,
and provide some necessary background.

\subsection{Notation and Setting}

 
\spara{Skills.} We consider a set $\skillset$ of skills with $\Card{\skillset}=\numskills$.
Skills can be any kind of qualification a worker can have or a task may
require, such as
\emph{video editing}, \emph{technical writing}, or \emph{project management}. 

\spara{Tasks.} We consider a set of 
$\numtasksTotal$ tasks (or jobs), $\taskset = \{\task{t}\St t=1,2,\dots, \numtasksTotal\}$, 
which arrive one-by-one in a streaming fashion;
$\task{t}$ is the $t$th task that arrives. Each task $\atask\in\taskset$ requires a
set of skills from $\skillset$, therefore, $\atask\subseteq \skillset$.
We use $\task{t}$ to refer to both
the task and the skills that it requires.


\spara{Workers.} 
Throughout we assume that we have a set 
$\workerset$ of $\numworkers$ workers: $\workerset = \{
\worker{r}\St r = 1,\dots,\numworkers \}$.
Every worker $r$ possesses a set of skills ($\worker{r}\subseteq\skillset$), and 
$\pool{\ell}$ denotes the subset of workers possessing
a given skill $\ell$: $\pool{\ell}=\{r\St \ell \in \worker{r}\}$.
Similarly to the tasks,
we use $\worker{r}$ to denote both the worker and his/her skills.
%
%

We partition the set of available workers $\workerset$ into the set of
workers who
are \emph{hired} at time $t$, denoted by $\hiredworkers^t$, and the set of workers
who are \emph{not hired}, denoted by $\nonhiredworkers^t$ (we sometimes
refer to these workers as \emph{\underline{f}reelancers}, and they can be \emph{outsourced} for $\task{t}$), so that
$\hiredworkers^t\cap\nonhiredworkers^t=\emptyset$ and
$\workerset=\hiredworkers^t\cup\nonhiredworkers^t$.

\spara{Coverage of tasks.} Whenever task $\task{t}\subseteq\skillset$
arrives, an algorithm has to assign one or more workers to it, i.e., a \emph{team}.
We say that $\task{t}$ can be \emph{completed} or \emph{covered} by
a team $\solution\subseteq\workerset$
if for every skill required by $\task{t}$, there exists at least one
worker in $\solution$ who possesses this skill:
$\task{t} \subseteq \cup_{\aworker\in\solution}\aworker$.
We assume that for every skill in the incoming task there is at least one worker possessing that skill, so all tasks can be covered.


\spara{Costs.}
Every worker $\worker{r}$ potentially can charge the
following nonnegative, worker-specific fees: 
\begin{inparaenum}
 \item an \emph{outsourcing fee} $\outsourcecost{r}$,
 \item a \emph{hiring fee} $\hiringcost{r}$, and
 \item a \emph{salary} $\salary{r}$.
\end{inparaenum}
Outsourcing fees $\outsourcecost{r}$ denote the payment required by a (non-hired) worker when a task is outsourced to him/her.
Note that $\outsourcecost{r}$ depends on the worker but does not depend on the task.
Hiring fees $\hiringcost{r}$ reflect all expenses associated to hiring
and firing a worker, such as signup bonuses and severance payments.
Given that any algorithm commits to pay the firing costs the moment in which it hires a worker, we follow a standard methodology used in online algorithms for caching~\cite{buchbinder2009design} and account for both hiring and firing costs when the worker is hired.
Once a worker $r$ is hired, s/he is paid a recurring salary $\salary{r}$,
which recurs for every step $t$ that the worker is hired. 
The above notation is summarized on Table~\ref{tbl:notation}.

\spara{Assumptions.}
To avoid making the model overly complicated, we assume that the salary periods are defined by the arriving tasks, this is, there is one task per salary period, and task completion takes one salary period.
%
%
A further assumption will be that $\salary{r} < \outsourcecost{r}$, as in
practice requesting a single task from an external worker involves extra costs~\cite{barthelemy2001hidden}, which are reduced when the worker is hired (or when an outsourcing arrangement for an external group of workers to perform a specific recurring task is done, which is different from the individual outsourcing we discuss here).
Finally, we assume $\outsourcecost{r} < \hiringcost{r}+\salary{r}$,
because otherwise workers would be hired and fired for every task.

%

\begin{table}[t]
\caption{Notation}
\label{tbl:notation}
\centering\small\begin{tabular}{cl}
\toprule
$\skillset$	& Set of skills, size $\numskills$ \\
$\taskset$	& Set of tasks, size $\numtasksTotal$ \\
$\numtasks$	& Number of tasks till current time \\
$\task{t}$ & The $t$'th task arriving\\
		& $\task{t}_\ell = 1$ if task $t$ requires skill $\ell$, 0 otherwise\\
$\workerset$	& Set of workers, size $\numworkers$. \\
		& $\worker{r}_\ell=1$ if worker $r$ possess skill $\ell$, 0 otherwise\\
$\pool{\ell}$	& Subset of workers possessing skill $\ell$ \\
\midrule
$\hiringcost{r}$	& Hiring fee, paid when worker $r$ is hired \\
$\outsourcecost{r}$& Outsourcing fee, paid every time $r$ performs a task \\
$\salary{r}$ & Salary paid to a hired worker $r$\\
\bottomrule
\end{tabular}
\end{table}

\subsection{Problem Definition}

We now define the problem that we study:

\begin{problem}[Team Formation with Outsourcing -- \problemsalary]\label{problem:theproblem}
There exists a set of skills \skillset. We have a pool of workers
\workerset, where each worker $\worker{r}\in\workerset$ is characterized by
\begin{inparaenum}[(1)]
a subset of skills $\worker{r}\subseteq\skillset$,
an outsourcing cost $\outsourcecost{r}\in\mathbb{R}_{\ge0}$,
a hiring cost $\hiringcost{r}\in\mathbb{R}_{\ge0}$, and
a salary cost $\salary{r}\in\mathbb{R}_{\ge0}$.
\end{inparaenum}
Given a set of tasks $\taskset=\{\task{1},\task{2},\dots,\task{\numtasksTotal}\}$, with
$\task{t}\subseteq\skillset$, which arrive in a streaming
fashion,
the goal is to design an algorithm that, when task $\task{t}$
arrives, decides which workers to hire (paying cost
$\hiringcost{r}+\salary{r}$), keep hired (paying cost \salary{r}), and
outsource (paying cost \outsourcecost{r}),
such that all the tasks are covered by the workers who are hired or outsourced and the total cost paid over all the tasks is minimized.
\end{problem}

{\problemsalary} is an online problem: $\taskset$ is revealed one task at a time.
Our goal 
is to guarantee that
for \emph{any} input stream $\taskset$ the total cost of our online algorithm, ${\mathit ALG}(\taskset)$, is at most a small factor greater than the total cost
of the optimal (offline) algorithm that knows $\taskset$ in advance, ${\mathit OPT}(\taskset)$.
%
%
This factor, $\max_{\taskset} {\mathit ALG}(\taskset)/{\mathit OPT}(\taskset)$, is called the \emph{competitive ratio} of the algorithm.


We solve the \problemsalary problem in Section~\ref{sec:salarynew}. Because neither the
algorithm nor its analysis are trivial, we introduce them gradually by first 
solving a simplified version of {\problemsalary}, which we describe
and solve in Section~\ref{sec:useforfree}.

\subsection{Background Problems}
\label{sec:background}

Two special cases of \problemsalary are \setcoverproblem and \skirentalproblem. 

{\bf\setcoverproblem}: {\bf the single-task, multiple-skill case.}
The set cover problem is an instance of our problem
when there is a single task $\atask\subseteq\skillset$ and
for each worker $\worker{r}$, $\hiringcost{r}=\infty$.
Then, as soon as the task \atask arrives, the algorithm
needs to cover all skills in $\atask$ by selecting a
set of workers $\solution\subseteq \workerset $ such that $\solution$
covers $\atask$ and 
$\sum_{r\in\solution}\outsourcecost{r}$ is minimized.
In this case, our problem can be solved using the greedy algorithm
for the set-cover problem (see~\cite[Chapter 2]{vazirani2013approximation}).

{\bf\skirentalproblem}: {\bf the single-skill, single-worker case.}
The ski rental problem is an instance of our problem
when the sequence of tasks $\taskset$ consists of a repetition
of the same single-skill task $\atask$ and the workforce $\workerset$
consists of a single worker $\worker{r}$ who possesses the same one skill,
and has $\salary{r}=0$ and some $\hiringcost{r}$, $\outsourcecost{r}$.
In this ski-rental 
version of our problem~\cite{manasse2008ski},
the question is the following:
without knowledge of
the total number of tasks that will
arrive, when should worker $\aworker$ be
hired so that the total cost paid to him/her in outsourcing plus hiring fees
is minimized?

A well-known algorithm for this problem is the following:
for every instance of $\task{t}$ that arrives outsource $\task{t}$ to worker
$\worker{r}$ as long as: $\sum_{t'=1}^t\outsourcecost{r}<\hiringcost{r}$.  Then, hire the worker when $\sum_{t'=1}^t\outsourcecost{r}\geq
\hiringcost{r}$. 
The above algorithm achieves a
competitive ratio of~$2$.


%% file: simple.tex
\newcommand{\poolUnhired}[1]{\ensuremath{P^{\nonhiredworkers}_{#1}}\xspace}
\newcommand{\poolUnhiredUseful}[1]{\ensuremath{P^{\nonhiredworkers}_{#1}}\xspace}

First, we solve a simplified version of the {\problemsalary} problem, where for every worker $\worker{r}$ the salary is equal to 0 ($\salary{r}=0$).
%
In this version of the problem, which we call {\problemuseforfree}, 
a hired worker $\worker{r}$ is paid a lump sum of $\hiringcost{r}$ the moment s/he is hired and this amount is assumed to cover all future work done by the worker.
Instead, when a worker $\worker{r}$ is outsourced, a payment of
$\outsourcecost{r}$ is done every time s/he performs a task.

\subsection{The {\naive} Algorithm}
\label{subsec:simple-naive}

A natural algorithm for solving the {\problemuseforfree} problem
combines ideas from {\setcoverproblem} and {\skirentalproblem} as follows:
first, it starts with no worker being hired and
each worker $\worker{r}$ is associated with a
variable $\delta_r$ initially set to~$0$.

For any $\numtasks\in\{1, \dots, \numtasksTotal\}$,
when task $\task{\numtasks}$ arrives, the algorithm
proceeds as follows: first, it identifies $\task{\numtasks}_\nonhiredworkers$
to be the set of skills of $\task{\numtasks}$ that cannot be covered
by already-hired workers.
Then, it covers the skills in $\task{\numtasks}_\nonhiredworkers$
using the greedy algorithm for \setcoverproblem. This way it finds
$\solution^\numtasks\subseteq\workerset$ such that
$\sum_{\worker{r}\in\solution^\numtasks}\outsourcecost{r}$ is minimized.
Finally, for each worker $\worker{r}\in\solution^\numtasks$, it updates
$\delta_r \leftarrow \delta_r+\outsourcecost{r}$. Worker $\worker{r}$ is hired when
$\delta_r\geq \hiringcost{r}$.
Clearly, since there are no salaries there is no motivation to fire a worker
once s/he is hired.

\spara{\naive has arbitrarily bad competitive ratio.}
Although our experiments (Section~\ref{sec:experiments})
demonstrate that the above algorithm, which we call {\naive}, performs
quite well in many practical cases, 
we can show that its competitive ratio can be arbitrarily bad.
For this, consider an example where $\workerset=\{\worker{1},\worker{2}\}$
and both workers have the same skill: $\worker{1}=\worker{2}=\{\askill\}$.
Further assume that $\outsourcecost{1}=1$, $\outsourcecost{2}=1+\epsilon$
and $\hiringcost{1}=M$, $\hiringcost{2}=2$, where $M$ is a large
value and $\epsilon$ a small one.
For a sequence of tasks $\task{1}=\task{2}=\ldots = \task{\numtasksTotal}=\{\askill\}$, it is clear that {\naive} will always outsource to $\worker{1}$ until
hiring him/her and will incur worst-case cost $2M$, whereas the optimal algorithm pays just~$\hiringcost{2}=2$. 

\subsection{A Primal--Dual Algorithm}
The above discussion illustrates that to obtain an algorithm with bounded
competitive ratio,
we need to take into
account both the outsourcing and hiring costs of all workers.
To do so, we 
deploy an online primal--dual scheme, which drives our algorithm design.

\spara{The integer and linear programs.}
The first step of the primal--dual approach, is to define an
integer formulation for the problem, for each step
$\numtasks\in\{1, \dots, \numtasksTotal\}$. We assume that the current task is
the $\numtasks$th task and we use the following variables:
\begin{compactitem}
\item $x_r=1$ if worker \worker{r} is hired when task $\task{T}$
arrives; otherwise $x_r=0$.
\item $f_{rt}=1$ if worker \worker{r} is outsourced for performing
task $\task{t}$; otherwise $f_{rt}=0$.
\end{compactitem}

Using this notation, {\problemuseforfree} can be formulated
as follows:

\begin{framed}
Linear program for \problemuseforfree:
\[
\min \sum_{r = 1}^{\numworkers}  \left(\hiringcost{r} x_r + \outsourcecost{r}\sum_{t=1}^\numtasks f_{rt}\right)
\hspace{2cm}\text{subject to:}
\]
\begin{flalign}\label{eq:tfo-lump-sum-primal}
\text{$\forall t = 1,\dots,\numtasks, \askill \in \task{t}:$}\qquad
&\sum_{\worker{r}\in \pool{\askill}} \left(x_r+f_{rt}\right)  \geq  1&
\end{flalign}
\begin{flalign*}
\text{$\forall t = 1,\dots,\numtasks, r=1,\dots,\numworkers$:}\qquad
&x_r, f_{rt}  \ge  0&
\end{flalign*}
\end{framed}
\medskip

The above, in addition to the integrality constraints $x_r, f_{rt}\in\mathbb{N}$, form the integer program for {\problemuseforfree}.
%
In this formulation, the objective function
sums over all workers the hiring costs (paid if the corresponding
worker has been hired by time $t$) and the outsourcing cost for the
tasks for which the worker has been outsourced. This is the total cost
of the solution until the current task~\task{\numtasks}.
Note that in this formulation of the problem there is no motivation for
a worker who is hired to be fired. Therefore, once $x_r$ is set to $1$,
it does not change its value to become $0$ again.

The first constraint \eqref{eq:tfo-lump-sum-primal} in the above program is the covering constraint: it simply enforces that
for every skill required for each task, there exists a hired or outsourced
worker who has this skill.
This guarantees that the team selected for
each task \task{t} covers \emph{all} the required skills.
The nonnegativity and the integrality constraints,
ensure that the solutions that we obtain from the integer-program
formulation can be transformed to a solution to our problem:
eventually, every variable will take the value 0 or~1.\footnote{A
solution in which some variables take values greater than 1, can be
transformed to another feasible solution with lower cost by setting
these variables to~1.}

To apply the online primal--dual approach, we first consider the linear
relaxation of the integer program, which simply drops the integrality
constraints $x_r, f_{rt}\in\mathbb{N}$. In a solution to this linear program
(LP) each variable takes values in
$[0,1]$. 
Given this LP, we can write its dual as follows:


\medskip

\begin{framed}
The dual of the linear program for \problemuseforfree:
\[
\max \sum_{t=1}^\numtasks \sum_{\askill\in\task{t}} u_{\askill t}
\hspace{2cm}\text{subject to:}
\]
\begin{flalign}\label{eq:skirental-buy}
\text{$\forall r=1,\dots,\numworkers $:}\qquad
\sum_{t=1}^\numtasks \sum_{\askill \in \task{t} \cap \worker{r}}
u_{\askill t}&\leq  \hiringcost{r}&
\end{flalign}
\begin{flalign}\label{eq:skirental-rent}
\text{$\forall t=1,\dots,\numtasks , r=1,\dots,\numworkers$:}\qquad
&\sum_{\askill\in \task{t} \cap \worker{r}} u_{\askill t}\leq\outsourcecost{r}&
\end{flalign}
\begin{flalign*}
\text{$\forall t = 1,\dots,\numtasks, \askill \in \task{t}$:}\qquad
&u_{\askill t}  \ge  0&
\end{flalign*}
\end{framed}
\medskip

Note that at every time $t \in\{1,\dots,\numtasks\}$
we have such a pair of primal--dual
formulations.
We are now going to use these
two formulations for designing and analyzing our algorithm.

\spara{The {\algouseforfree} algorithm:}
%
Next, we present the {\algouseforfree} algorithm, which is designed and analyzed using the primal and the dual linear programs. 
%
We assume that task $\task{\numtasks}$, 
$\numtasks\in\{1, \dots, \numtasksTotal\}$, has just arrived and the
algorithm must act before task $\task{\numtasks+1}$ arrives (or the stream
finishes if $\numtasks=\numtasksTotal$).

All the variables used in our algorithm are initialized to 0 before the
arrival of the first task.
When task $\task{\numtasks}$ arrives the following steps are done:

\begin{enumerate-algo}

\item Let $\nonhiredworkers^\numtasks$ and $\hiredworkers^\numtasks$ represent the
workers who are \emph{not hired} and \emph{hired}, respectively, at the time that
$\task{\numtasks}$ arrives. Clearly, when the first task arrives ($\numtasks=1$), 
$\nonhiredworkers^\numtasks=\workerset$ and
$\hiredworkers^\numtasks=\emptyset$.
For $\numtasks>1$, the values of $\hiredworkers^\numtasks$ and
$\nonhiredworkers^\numtasks$ are updated in the last step (step 10) of
the previous round.


\item Let
$\task{\numtasks}_\hiredworkers=\task{\numtasks}\cap\cup_{\worker{r}\in\hiredworkers^\numtasks}\worker{r}$
be the skills from
$\task{\numtasks}$ that are covered by already-hired workers and
$\task{\numtasks}_\nonhiredworkers = \task{\numtasks}\setminus
\task{\numtasks}_\hiredworkers$.
  
\item For every skill $\askill\in\task{\numtasks}_\nonhiredworkers$ let
$\poolUnhired{\askill}=\pool{\askill}\cap\nonhiredworkers^\numtasks$
be the set of workers in $\nonhiredworkers^\numtasks$ such
that every worker in $\poolUnhired{\askill}$ has skill $\askill$.
Also let 
\(\poolUnhiredUseful{\numtasks}=\cup_{\askill\in\task{\numtasks}_\nonhiredworkers}\poolUnhired{\askill}\)
be the set of unhired workers who possess at least one skill that is required and not covered by already-hired workers.



\item \textbf{for each}
$\worker{r}\in\poolUnhired{\numtasks}$:\hspace{.4cm}set
$\tilde{x}'_r\leftarrow \tilde{x}_r$.

\item \textbf{for each} skill $\askill\in\task{\numtasks}_\nonhiredworkers$:

\hspace{0cm} \textbf{while} $\sum_{\worker{r}\in \pool{\askill}}\left(\tilde{x}_r+\tilde{f}_{r\numtasks}\right)<1$:


\hspace{.4cm}\textbf{for each} $\worker{r}\in \pool{\askill}$: $\tilde{x}_r\leftarrow \tilde{x}_r\left(1+\frac{1}{\hiringcost{r}}\right)+\frac{1}{\numworkers\hiringcost{r}}$

\hspace{.4cm}\textbf{for each} $\worker{r}\in \pool{\askill}$: $\tilde{f}_{r\numtasks}\leftarrow \tilde{f}_{r\numtasks}\left(1+\frac{1}{\outsourcecost{r}}\right)+\frac{1}{\numworkers\outsourcecost{r}}$

\item \textbf{for each}
$\worker{r}\in\poolUnhired{\numtasks}$:\hspace{.4cm} set $\Delta \tilde{x}_r\leftarrow \tilde{x}_r - \tilde{x}'_r$.

\item Set $\hiredworkers'\leftarrow \emptyset$. 


\item \textbf{repeat} $\rho_1$ times:
    
      \hspace{0.0cm}\textbf{for each}
$\worker{r}\in\poolUnhiredUseful{\numtasks}$
    
        \hspace{0.4cm}with probability $\Delta \tilde{x}_r$:
    
          \hspace{0.8cm}hire worker $\worker{r}$ (set
    	  $x_r\leftarrow1$, $\hiredworkers'\leftarrow\hiredworkers'\cup\{r\}$)

        \hspace{0.4cm}with probability $\tilde{f}_{r\numtasks}$:

          \hspace{0.8cm}outsource worker $\worker{r}$ (set $f_{r\numtasks}\leftarrow1$)

%
%
%


\item \textbf{for each} skill $\askill\in\task{\numtasks}_\nonhiredworkers$:

    \hspace{.0cm}\textbf{if} skill $\askill$ is not covered:

    \hspace{.4cm}hire worker $\worker{r}\in\poolUnhired{\askill}$ with minimum cost $\hiringcost{r}$

    \hspace{.4cm}(set $x_r\leftarrow1$, $\hiredworkers'\leftarrow\hiredworkers'\cup\{r\}$)
%
        
        
\item $\hiredworkers^{\numtasks+1}\leftarrow\hiredworkers^{\numtasks}\cup\hiredworkers',\quad
\nonhiredworkers^{\numtasks+1}\leftarrow\workerset\setminus\hiredworkers^{\numtasks+1}$.
\end{enumerate-algo}

For $\numtasks=1$, the {\algouseforfree} starts with no worker being hired.
Intuitively, as tasks arrive, the algorithm tries to gauge two
quantities:
(1) the usefulness of every
worker for the task at hand $\task{\numtasks}$ and (2) the overall
usefulness of every worker for tasks $\task{1},\ldots,\task{T}$. This is done
in step~5, via variables $\tilde{f}_{r\numtasks}$ (for (1)) and
$\tilde{x}_r$ (for (2)). In particular,
the more useful \worker{r} proves over time,
the larger the value~$\tilde{x}_r$.
Subsequently, in step~8 every worker is outsourced or hired based on
the increase in the values of  $\tilde{f}_{r\numtasks}$ and $\tilde{x}_r$
observed in step~5.
Finally, for every skill that remains uncovered after
step~8 (which is randomized), {\algouseforfree} hires worker $\worker{r}$ with
the minimum $\hiringcost{r}$ that covers the skill.
Note that the increase of the variables $u_{\askill\numtasks}$ in step~5 is
not required for solving the {\problemuseforfree},
but it is used in our analysis and thus we leave it in the description above.

Our analysis requires to set the value of $\rho_1$ in step~8
to
\(\rho_1=\ln\numskills+\ln\maxhiringcost{r}\),
where
$\maxhiringcost{r}=\max_{\worker{r}\in\workerset}\hiringcost{r}$.

Although one may think that an additive update of
variables in 
step~5 would seem more natural, such an update would introduce a
$\Theta(\numskills)$
factor in the competitive ratio. On the other hand, the multiplicative update
that we adopt, has the property that the more a worker \worker{r} is required over
time the higher the increase of the corresponding variable
$\tilde{x}_r$. This fact, leads us to Theorem~\ref{theorem:useforfree} below.

\spara{Analysis.}
We have the following result for \algouseforfree.

\begin{theorem}\label{theorem:useforfree}
  \algouseforfree is an
$O(\log\numworkers(\log\numskills+\log\maxhiringcost{r}))$- competitive algorithm
 for  the \problemuseforfree problem, where
 $\maxhiringcost{r}=\max_{\worker{r}\in\workerset}\hiringcost{r}$.
\end{theorem}

\spara{Running time.} The running time of {\algouseforfree} per task is dominated by the
execution of steps 5 and~8.
For step~5, using binary search, the algorithm can determine in $O(\log\maxhiringcost{r})$ steps the 
minimum increase of
$\tilde{x}_r$ and 
$\tilde{f}_{r\numtasks}$
that makes false the condition of the 
\texttt{while} loop for at least one uncovered skill~$\askill$.
Therefore, the running time of step~5
is $O\!\left(\Card{\task{\numtasks}}\numworkers\log\maxhiringcost{r}\right)$.
Step~8, using a hash table to store hired workers, can be executed in expected time
$O(\rho_1\numworkers)=O\left(\numworkers(\log\numskills+\log\maxhiringcost{r})\right)$.
Therefore, the expected time required for processing task $\task{T}$ is
$O\!\left(\numworkers\left(\log\maxhiringcost{r}\Card{\task{\numtasks}}+\log\numskills\right)\right)$.

%% file: salarynew.tex
In this section, we provide an algorithm for the general
version of {\problemsalary} (Problem~\ref{problem:theproblem}).
In contrast with {\problemuseforfree}, now after hiring a worker we must
pay a salary $\salary{r}\geq 0$, complicating the problem significantly
as it may now be cost-effective to fire workers.

\spara{The integer and linear programs for \problemsalary.} 
Given that workers can be hired, then fired and potentially hired again, and so on, we introduce in this new LP the notion of intervals.
These intervals are used to model periods in which workers are hired $\intervals=\left\{\left\{t_a,t_b\right\}\mid t_a,t_b\in\mathbb{N}, t_a\le t_b\right\}$.
Intuitively, an interval is a subset of time steps during which
an algorithm decides to hire a given worker.
The new LP, (omitted) 
uses the following variables:
\begin{compactitem}
\item $x(r,I)$ with $I\in \intervals$: $x(r,I)=1$ if worker $\worker{r}$ is hired during the entire interval $I$;
otherwise $x(r,I)=0$.
\item $f_{rt}$: $f_{rt}=1$ iff worker $\worker{r}$ is outsourced for performing
$\task{t}$.
\end{compactitem}

It turns out that it is hard to design an approximation algorithm with
proven guarantees using
this program, mostly because it is hard to keep track of the costs being paid
for every worker when the intervals of him/her being hired, outsourced, or idle are of variable length. 
Therefore, we resort to a different 
overall \emph{strategy:}
First, we define the {\problemalternative} problem, in which the solutions are restricted such that every worker
is hired for fixed-length (worker-specific) intervals (Section~\ref{subsec:problem-alternative}).
Then, we design an algorithm for
{\problemalternative} with
good competitive ratio (Section~\ref{subsec:alt-tfo-algorithm}).
Finally, we prove that a solution to {\problemalternative}
can be transformed to a solution for \problemsalary, and
that any solution of {\problemsalary} 
can be transformed to a feasible solution
of {\problemalternative} that is a factor of at most $3$ times higher (Section~\ref{subsec:alt-to-salary}), obtaining an approximation algorithm for
\problemsalary.

\subsection{The {\problemalternative}
  Problem}\label{subsec:problem-alternative}

The difference between {\problemalternative} and {\problemsalary}
is that we restrict the solutions of the former to have 
a specific structure; 
whenever worker $\worker{r}$ is hired s/he is
then fired after $\hireintervals_r\Def \Ceil{\hiringcost{r}/\salary{r}}$ time
units---independently of whether s/he is used or not in tasks within these
$\hireintervals_r$ time units.

In this case, every worker $\worker{r}$ is associated with a new hiring cost
$\hiringcostwithsalary{r}$, which is the summation of
his/her original hiring cost $\hiringcost{r}$ plus the salaries paid
to him/her for the $\hireintervals_r$ time units he is hired.
Thus, the total hiring cost and salary for an entire interval is
$\hiringcost{r}+\hireintervals_r\cdot\salary{r}\le
\hiringcost{r}+\left(\frac{\hiringcost{r}}{\salary{r}}+1\right)\cdot\salary{r}
\leq 3 \hiringcost{r}.$
We will use
$\hiringcostwithsalary{r}\Def 3 \cdot\hiringcost{r}$.

We can now write the LP for {\problemalternative}. 
In addition to the notation we discussed in the previous paragraph, we use 
$I^t\in \intervals$ to denote the 
interval that starts at time $t$.
Worker $\worker{r}$ has $x(r,I) = 1$ iff s/he is hired during
the \emph{entire} interval $I$.
All intervals $I$ for which $x(r,I)=1$ are of fixed length
$\hireintervals_r$. 

\begin{framed}
Linear program for {\problemalternative}:
\[
\min \sum_{r = 1}^{\numworkers} \left[ \sum_{I\in \intervals}  \hiringcostwithsalary{r} x(r,I)
 + \sum_{t=1}^{\numtasks} \outsourcecost{r}  f_{rt} \right]
\hspace{2cm}\text{subject to:}
\]
\begin{flalign}\label{eq:cover-skills-alt}
\text{$\forall t=1\ldots\numtasks, \askill\in\task{t}:$}\quad
&\sum_{\worker{r}\in P_\askill}\left(f_{rt} + \sum_{I\in \intervals: t\in I} x(r,I) \right) \geq 1&
\end{flalign}
\begin{flalign*}
\text{$\forall t=1\ldots\numtasks, r=1\ldots\numworkers,I\in\intervals:$}\quad
&x(r,I), f_{rt}\ge0&
\end{flalign*}
\end{framed}
\medskip

\subsection{Solving the {\problemalternative} Problem}\label{subsec:alt-tfo-algorithm}
In this section, we design and analyze an algorithm for the {\problemalternative} problem.
The similarity between the LPs for {\problemalternative} and {\problemuseforfree} (Section~\ref{sec:useforfree})
translates into a similarity in the algorithms (and their analysis) of
the two problems. The key difference now is that we need to take care of the firings.


Our algorithm for {\problemalternative} differs from the algorithm for {\problemuseforfree} in steps 1, 5, 8, and 9, which are changed 
as follows:

\begin{enumerate-algo}
\item[1'.]
Let $\nonhiredworkers^\numtasks$ and $\hiredworkers^\numtasks$ represent the
workers who are not hired and hired, respectively, at the time that
$\task{\numtasks}$ arrives. Clearly, when the first task arrives ($\numtasks=1$), then
$\nonhiredworkers^\numtasks=\workerset$ and
$\hiredworkers^\numtasks=\emptyset$.
For $\numtasks>1$, the values of $\hiredworkers^\numtasks$ and
$\nonhiredworkers^\numtasks$ are updated in the last step (step 10) of
the previous round and then we remove workers whose hiring interval
finished in the previous step:

\hspace{.4cm}$\mathcal{F}'\leftarrow\{r\in\workerset\St x(r,I^{\numtasks-\hireintervals_r}) = 1\}$

\hspace{.4cm}$\hiredworkers^{\numtasks}\leftarrow\hiredworkers^{\numtasks}\setminus\mathcal{F}',
\quad\nonhiredworkers^{\numtasks}\leftarrow\workerset\setminus\hiredworkers^{\numtasks}$

\hspace{.4cm}\textbf{for each} $\worker{r}\in\mathcal{F}'$: \hspace{.4cm}set $\tilde{x}_r\leftarrow0$

\item[5'.]
\textbf{for each} skill $\askill\in\task{\numtasks}_\nonhiredworkers$:

\hspace{.0cm} \textbf{while} $\sum_{r\in \pool{\askill}}\left(\tilde{x}_r+\tilde{f}_{r\numtasks}\right)<1$:


\hspace{.4cm}\textbf{for each} $r\in \pool{\askill}$:\hspace{0.4cm}$\tilde{x}_r\leftarrow \tilde{x}_r\left(1+\frac{1}{\hiringcostwithsalary{r}}\right)+
\frac{1}{\numworkers\hiringcostwithsalary{r}}$

\hspace{.4cm}\textbf{for each} $r\in \pool{\askill}$:\hspace{0.4cm}$\tilde{f}_{r\numtasks}\leftarrow \tilde{f}_{r\numtasks}\left(1+\frac{1}{\outsourcecost{r}}\right)+
\frac{1}{\numworkers\outsourcecost{r}}$

\item[8'.] \textbf{repeat} $\rho_2(\numtasks)$ times:

      \hspace{0.0cm}\textbf{for each}
$r\in\poolUnhiredUseful{\numtasks}$

        \hspace{0.4cm}with probability $\Delta \tilde{x}_r$:

          \mbox{\hspace{0.8cm}hire worker $\worker{r}$ (set
          $x(r, I^\numtasks)\leftarrow1$,
          $\hiredworkers'\leftarrow\hiredworkers'\cup\{r\}$)}

        \hspace{0.4cm}with probability $\tilde{f}_{r\numtasks}$:

          \hspace{0.8cm}outsource worker $\worker{r}$ (set $f_{r\numtasks}\leftarrow1$)

\item[9'.] \textbf{for each} skill $\askill\in\task{\numtasks}_\nonhiredworkers$:

    \hspace{.0cm}\textbf{if} skill $\askill$ is not covered:

    \hspace{.4cm}\mbox{outsource worker $\worker{r}$, $r\in\poolUnhired{\askill}$,
    with minimum cost $\outsourcecost{r}$}

    \hspace{.4cm}(set $f_{r\numtasks}\leftarrow1$)

\end{enumerate-algo}

Our analysis requires to set $\rho_2(\numtasks)=\ln\numskills+\ln\maxoutsourcecost{r} +2\ln\numtasks,$
where $\maxoutsourcecost{r}=\max_{\worker{r}\in\workerset}\outsourcecost{r}$.



\spara{Analysis of {\altalgo}.}
Algorithm \altalgo gives a solution with proven theoretical
guarantees for \problemalternative.
As before, the multiplicative update is needed to obtain this competitive ratio.
We have the following theorem (proof omitted due to space constraints).



\begin{theorem}\label{thm:salary}
\altalgo is an $O(\log\numworkers(\log
      \numskills+\log\maxoutsourcecost{r}+\log\numtasksTotal))$-competitive algorithm
 for  the \problemalternative problem.
\end{theorem}


\subsection{Solving {\problemsalary} Using {\problemalternative}}
\label{subsec:alt-to-salary}
Note that any solution output by {\altalgo}
can be transformed into a feasible solution to 
the original {\problemsalary}
problem by setting $g_{rt} \leftarrow 1$ for each $r, t\in I$ for which
$x(r,I)=1$, and $g_{rt} \leftarrow 0$ otherwise.
We call the algorithm that runs {\altalgo} and subsequently does this transformation a its final step,
the {\algowithsalary} algorithm.

The question is whether {\algowithsalary} provides a solution with bounded competitive
ratio for the {\problemsalary} problem.
We answer this question affirmatively by showing
\begin{inparaenum}[(1)]
\item that the solution of \algowithsalary for the \problemsalary problem is feasible and has a cost bounded by the cost of \altalgo for the \problemalternative problem, and
 \item that any
solution for the {\problemsalary} problem can be turned into a feasible solution
to the {\problemalternative} problem at the expense of a small loss in the approximation factor.
\end{inparaenum}
These two suffice to prove that the solution produced
by \algowithsalary is a good solution for the \problemsalary problem.
We have the following result (proof omitted due to space constraints):

\begin{theorem}\label{thm:salaryfinal}
\algowithsalary is an $O(\log\numworkers(\log
      \numskills+\log\maxoutsourcecost{r}+\log\numtasksTotal))$-competitive algorithm
 for  the \problemsalary problem.
\end{theorem}

\spara{Running time.}
Similarly to Section~\ref{sec:useforfree}, the expected time required to
process task $\task{T}$ is
$O\!\left(\numworkers\left(\log\maxhiringcost{r}\Card{\task{\numtasks}}+\log\numskills+\log\numtasks\right)\right)$.

\spara{Lower bound.}
Note that there is little hope for significant improvement of our theoretical results.
In particular, Alon et al.~\cite{alon2009online} have proven a lower
bound of
$\Omega\!\left(\frac{\log\numworkers\log\numskills}{\log\log\numworkers+\log\log\numskills}\right)$
on the competitiveness of any deterministic algorithm for the \emph{unweighted
online set cover problem}.
The unweighted online set cover problem, is a special case of
\problemsalary (and of \problemuseforfree) where for each worker $\worker{r}$ we have
$\hiringcost{r}=\outsourcecost{r}=1$, $\salary{r}=0$, and 
for each task $\task{\numtasks}$ we have
$\task{\numtasks-1}\cup\{\askill\}$,
for some skill $\askill\in\skillset\setminus \task{\numtasks-1}$ (with
$\task{0}=\emptyset$).

%

\subsection{The \naivetfo}
\label{concept:naivetfo}
Similarly to \problemuseforfree, we
also consider the heuristic \naivetfo, which is a generalization
of \naive, for general
values of $\salary{r}$. Specifically,
the difference is that worker $\worker{r}$ is hired when
$\delta_r \ge C_r + \eta_r\cdot\sigma_r$, and is fired after
$\hireintervals_r$ tasks (see Sections~\ref{subsec:simple-naive} and \ref{subsec:problem-alternative} for definitions of $\delta_r$ and $\eta_r$).
Note that theoretically \naivetfo may perform arbitrarily bad: the
example of Section~\ref{subsec:simple-naive} holds for \naivetfo for
small $\salary{r}$. Yet, in Section~\ref{sec:experiments} we observe that
even though it does not offer the theoretical guarantees of
\algowithsalary, it performs well in practice.

\subsection{The \metaalgo algorithm}
\label{concept:metaalgo}


As we will see in Section~\ref{sec:experiments}, although \algowithsalary
gives theoretical guarantees for the worst-case performance, in practice
some of our other algorithms for the \problemsalary problem may perform better under some input
parameters. Given the low running time of all our solution approaches to
\problemsalary, we implemented the
\metaalgo algorithm.
This algorithm runs in parallel all the presented methods for
solving the \problemsalary problem (\algowithsalary, \naivetfo,
\algoAlwaysOutsource and \algoAlwaysHire), and selects at each time
the current minimum-cost algorithm to apply to solve the current task,
switching between algorithms when it is advantageous.
The asymptotic worst-case results hold for the 
\metaalgo algorithm as well. Furthermore, our experiments (see
Section~\ref{sec:experiments}) show that
it is beneficial to change the hiring policy even if we pay switching
costs.

%% file: experiments.tex
%
Our experiments seek to compare the total cost that would be incurred by companies using different algorithms to assign workers to a stream of incoming tasks.
We use synthetic datasets representing possible workloads, built using actual task requirements and worker skills from three large online marketplaces.
Synthetic data, while having the limitation of not reflecting the particular conditions of a specific company, allows us to evaluate the effectiveness of our algorithms under a broad range of conditions.
Section~\ref{subsec:datasets} introduces our datasets,
Section~\ref{subsec:results-lumpsum} presents results on the \problemuseforfree problem, and
Section~\ref{subsec:results-tfo} on the \problemsalary problem.

\subsection{Datasets}\label{subsec:datasets}

\begin{table}[t]
\caption{
Characteristics of the three source datasets used to generate workloads for our experiments. Numbers in italics correspond to tasks generated for the Upwork dataset, as explained in Section~\ref{subsec:datasets}.
}
\label{tbl:data-summary}
\centering\small\begin{tabular}{lrrr}\toprule
Dataset & \multicolumn{1}{c}{UpWork} & \multicolumn{1}{c}{Freelancer} & \multicolumn{1}{c}{Guru} \\\midrule
Skills ($m$)  &		2,335 & 175	& 1,639	 \\
Workers ($n$) &		18,000 & 1,211	& 6,119	 \\
Tasks ($T$)	  &		\textit{50,000} & 992	& 3,194	 \\
... distinct  &		\textit{50,000} & 600	& 2,939	 \\
... avg. similarity (Jaccard) 	& \textit{0.095}	& 0.045	& 0.018	 \\ \midrule
Average Skills/worker & 6.29     & 1.45      & 13.07 \\
Average Skills/task & \textit{41.88}    & 2.86      & 5.24 \\
\bottomrule
\end{tabular}
\end{table}

\begin{figure}[t]
\subfigure[\problemuseforfree UpWork]{\includegraphics[width=.23\textwidth]{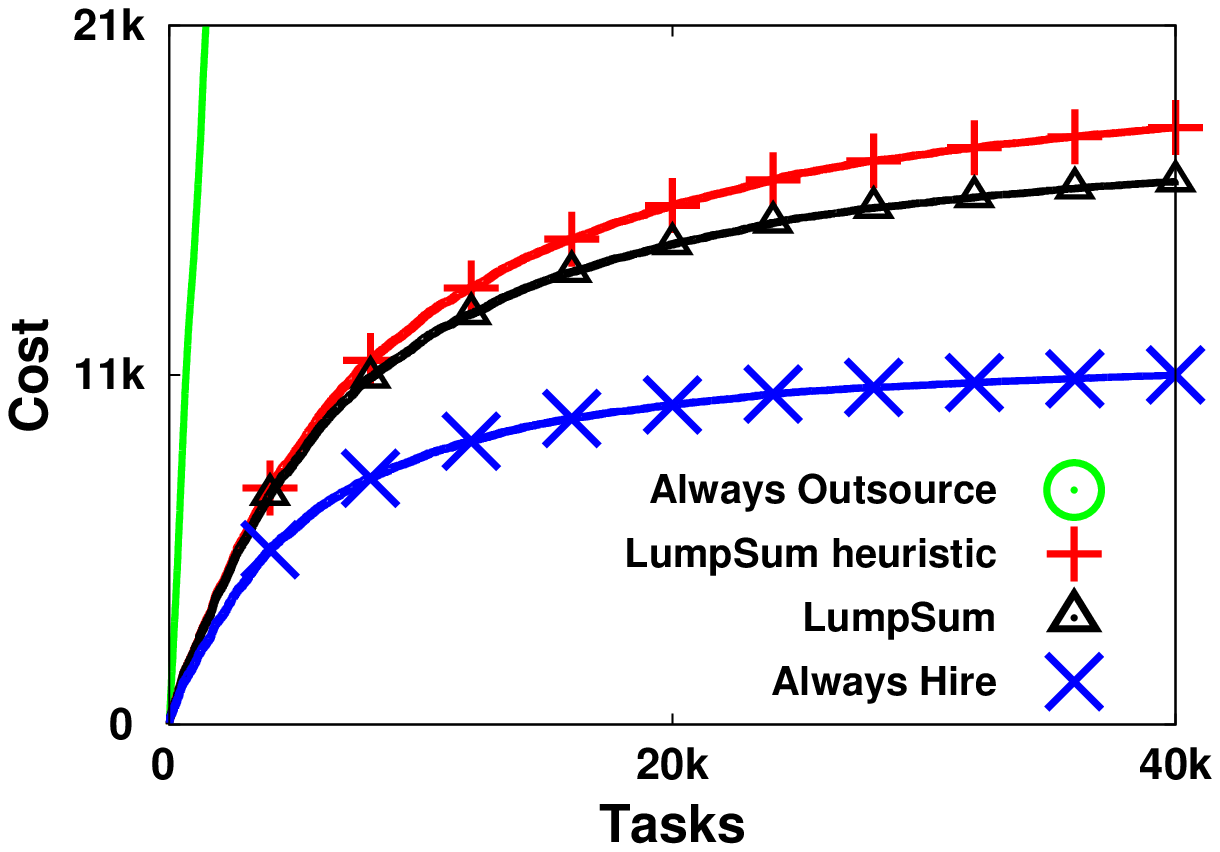}}
\subfigure[\problemsalary UpWork]{\includegraphics[width=.23\textwidth]{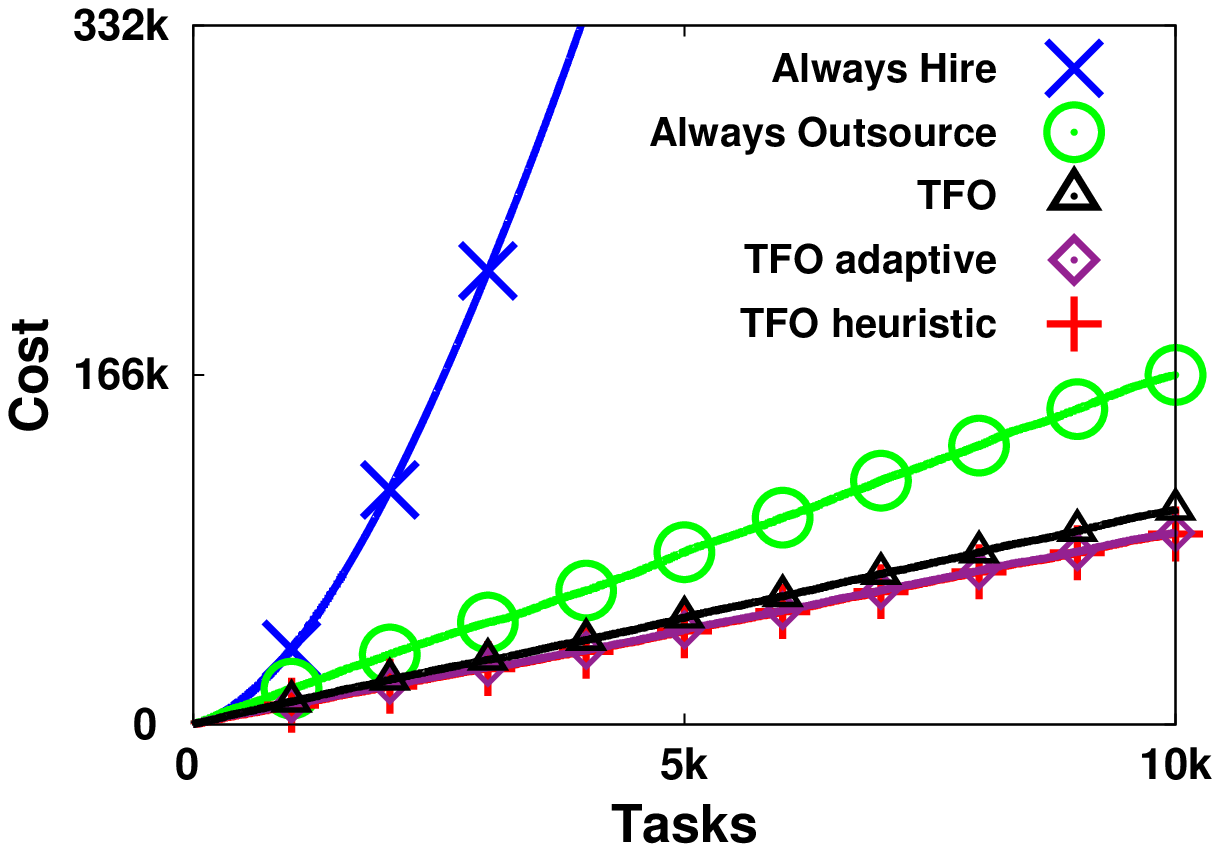}}
\subfigure[\problemuseforfree Freelancer]{\includegraphics[width=.23\textwidth]{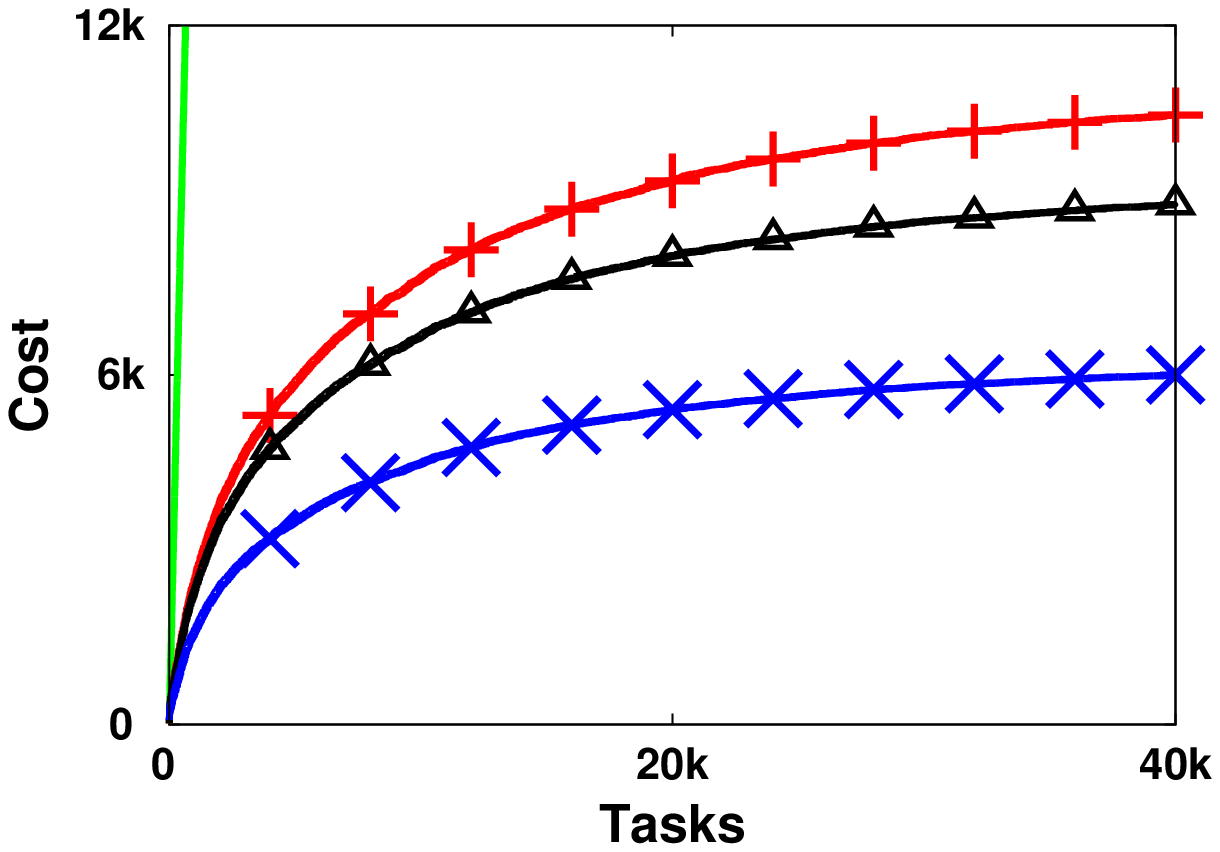}}
\subfigure[\problemsalary Freelancer]{\includegraphics[width=.23\textwidth]{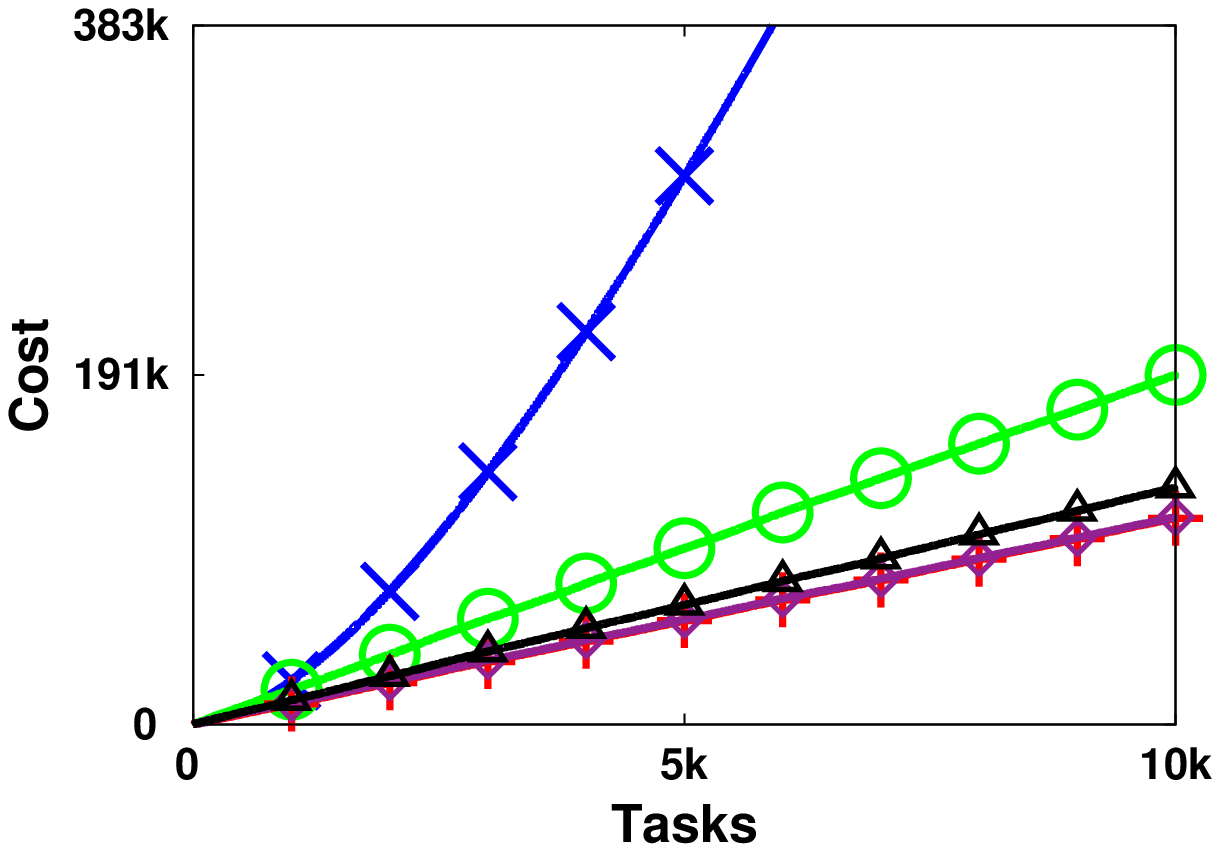}}
\subfigure[\problemuseforfree Guru]{\includegraphics[width=.23\textwidth]{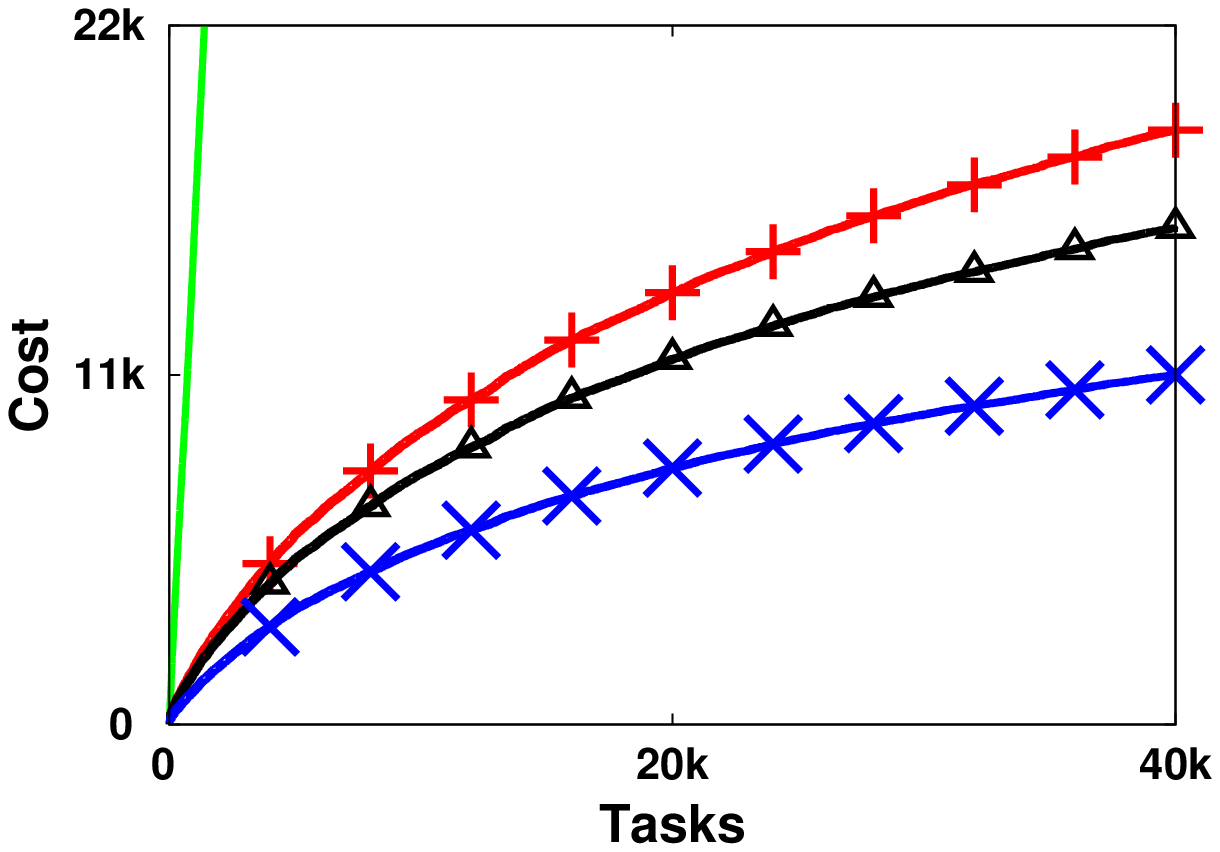}}
\subfigure[\problemsalary Guru]{\includegraphics[width=.23\textwidth]{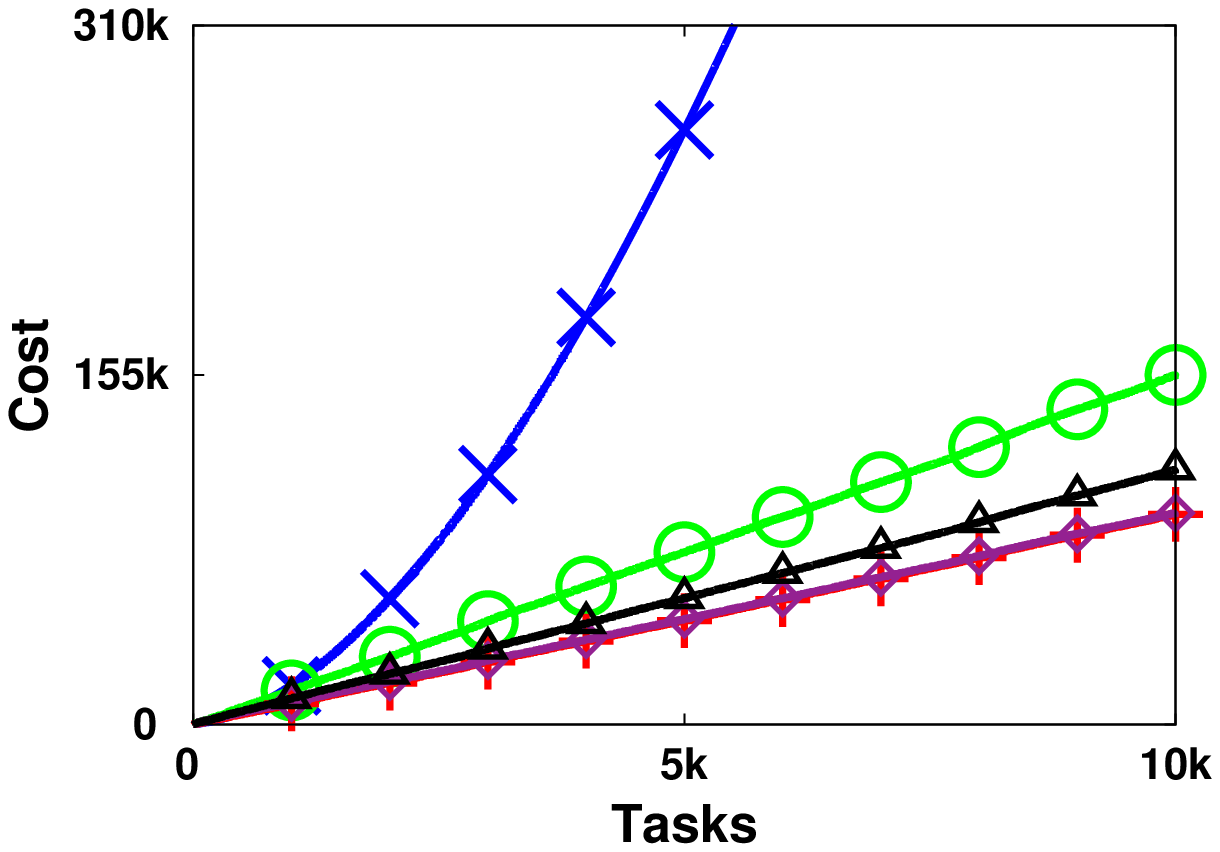}}
\caption{
Experimental comparison of algorithms showing total cost due to outsourcing, hiring, and paying salaries as a function of the number of tasks in the input, averaged over 100 workloads generated with $p=100$.
Left:
Algorithms for problem \problemuseforfree.
As expected, \algoAlwaysHire has the smallest cost if the number of tasks is large, however an online algorithm does not know the number of tasks. 
Our online algorithm and its heuristic version (\naive) show a cost that does not exceed twice of that of \algoAlwaysHire. In contrast, \algoAlwaysOutsource has cost proportional to the number of tasks.
Parameters $\hiringcost{r} = 4\outsourcecost{r}$ and $T=40\mathrm{K}$.
%
%
Right:
Algorithms for problem \problemsalary.
%
%
Our online algorithm, its heuristic version (\naivetfo) and the \metaalgo have smaller cost than \algoAlwaysOutsource and \algoAlwaysHire. The latter diverges rapidly due to salary costs. Parameters $\hiringcost{r} = 4\outsourcecost{r}$, $\salary{r} = 0.1\outsourcecost{r}$ and $T=10\mathrm{K}$.
%
\label{fig:results-problem-no-salary}
\label{fig:results-problem-salary}
}
\end{figure}

We start by introducing our datasets and discussing our choice of cost
parameters for experimentation.

\spara{Source datasets.}
To create a large pool of tasks from which to sample workloads, we use datasets obtained from three large online marketplaces for outsourcing:
UpWork, Freelancer and Guru (the authors are not associated with any of these services). 
All three are in the top-30 of traffic in their category (``consulting marketplaces'') according to data from Alexa (Feb. 2018),  
indeed, Freelancer and Guru are respectively number 1 and number 3.
%
%
General statistics of these datasets are shown on Table~\ref{tbl:data-summary}.


\spara{Worker skills.} The input data that we obtained contain anonymized profiles
for people registered as freelancers in these marketplaces. These
include their self-declared sets of skills, as well as the average rate
that they charge for their services.
%
There is a large variation in the number of skills per worker among datasets, as can be seen in
Table~\ref{tbl:data-summary}.
Data have been cleaned to remove skills that were not possessed by any worker
and skills that were never required by any task. The numbers in
Table~\ref{tbl:data-summary} refer to the clean datasets.

\spara{Tasks.} 
For both Freelancer and Guru we have access to a large sample of tasks
commissioned by buyers in the marketplace; they are included as tasks on
Table~\ref{tbl:data-summary}.
They correspond to actual tasks brought to these marketplaces by 
actual users.
%
%
These samples are anonymized: we do not know the name of the company commissioning them, and there are no timestamps in this data.
In the case of Upwork, we generate synthetic tasks following a data-generation procedure used in previous work~\cite{anagnostopoulos12online}: 
we remove a small number of workers (10\%), who are excluded from the pool of workers in the dataset,
and then repeatedly sample subsets of them to create
tasks, by interpreting the union of their skills as task requirements. 

\spara{Workloads.}
Marketplaces for online work cover a broad range of tasks from graphic design and web development to accounting, administrative assistance, and legal consulting. Except for huge conglomerates, most firms will not outsource work across all categories at the same time.
The workload-generation process that we use has a single parameter $p$, which we call the \emph{coherence parameter} of the workload, and works as
follows.
First, we start with a random task, which we select as
pivot. To select the next task, with probability $1/p$ we select a random task from the pool of distinct tasks in the dataset and make this task the new pivot, and with probability $1 - 1/p$ we select another task with Jaccard similarity at least $0.5$ to the pivot,
The expected length of a sequence of ``similar'' tasks is $p$.
Each workload stream that we create has 10K tasks. We also experimented with streams of up to $100$K tasks, but we observed that 10K tasks suffice to expose the trends of the algorithms that we compare.
%
We believe that in general a large value of $p$ is realistic for a company, as customers would probably procure from it services exhibiting a certain coherence; we also evaluate our algorithms for a broad range of values for~$p$.

For each dataset and for each coherence parameter that we use, we generated 100 workload streams; the costs that we report in our experiments are averages over these 100 workloads.

\spara{Cost parameters.}
We have data about the rates charged by workers in each marketplace,
which we directly interpret as their outsourcing costs $\outsourcecost{r}$.
However, we do not have their hiring or salary costs, so we experiment with different values for these costs.

For \emph{hiring costs}, which are characterized by $\hiringcost{r} >
\outsourcecost{r}$, we assume they are a multiplicative factor larger than the hiring cost, $\hiringcost{r} = \alpha_r \outsourcecost{r}$.
We performed extensive experiments in which $\hiringcost{r}$ varied between $1 \outsourcecost{r}$ and $30 \outsourcecost{r}$, either as a fixed value, or setting $\alpha_r$ to be a random variable distributed uniformly in a small range.
For \emph{salary costs}, 
we assume that they are a fraction of outsource costs, experimenting with values from $\salary{r} = \outsourcecost{r}/100$ to $\salary{r} = \outsourcecost{r}/4$.
%
%
Salaries $\salary{r}$ are smaller than outsourcing costs $\outsourcecost{r}$ because the latter includes many costs in which a company incurs when outsourcing~\cite{barthelemy2001hidden}, including:
\begin{inparaenum}[(i)]
\item outside-hired consultants are usually more highly paid per 
hour/day than regular employees for a company, 
\item there are transaction costs involved in locating and contracting 
and outsourced worker that do not exist for regular employees, and
\item there are communication and management costs of handling someone external to a company.
 \end{inparaenum}


%
%
%


\subsection{Experiments with \problemuseforfree}\label{subsec:results-lumpsum}

\spara{Baselines.}
We consider two baselines.
The \algoAlwaysHire baseline solves the \setcoverproblem problem for
finding a low-cost set of workers that cover the task's uncovered skills
and hires them.
%
%
The \algoAlwaysOutsource baseline never hires, instead it outsources to
workers that cover the required skills for the task, by solving a
\setcoverproblem problem instance.
%

\spara{Results.}
Figure~\ref{fig:results-problem-no-salary} (Left) summarizes our results for
\problemuseforfree for workloads generated with the UpWork, Freelancer and Guru datasets, depicting total cost as a function of the number of tasks.
%

We observe that under all these workloads the algorithms behave similarly.
\algoAlwaysOutsource has cost proportional to the number of tasks and is
not competitive, its cost is mostly outside the range of
Figure~\ref{fig:results-problem-no-salary} (Left).
As expected, \algoAlwaysHire performs the best in the long run,
because if the number of tasks is large, hiring is a dominant
strategy; however the online algorithm does not know the number of
tasks.
%
%
Experimentally, the \algouseforfree algorithm has a cost that does not exceed that of \algoAlwaysHire by more than a factor of 2, across all the scenarios that we tested.
We note that for short sequences {\algouseforfree} 
has lower cost; this difference in the cost can sometimes be an order of magnitude smaller (plots omitted for brevity). 
We also note that although \naive can, theoretically, perform arbitrarily bad, in our experiments it performs
quite well---although worse than the theoretically justified \algouseforfree.

\spara{Variations (plots omitted for brevity).}
Figure~\ref{fig:results-problem-no-salary} (Left) is obtained with $\hiringcost{r} = 4 \outsourcecost{r}$.
We do not observe dramatic variations in the results when varying this parameter in the studied range ($1\outsourcecost{r}$ through $30\outsourcecost{r}$): \algouseforfree has a smaller cost than \algoAlwaysOutsource.
In general, higher hiring costs mean the number of tasks required before hiring a worker is larger, the costs of \algouseforfree and \algoAlwaysHire are higher, and the advantage for \algouseforfree over \algoAlwaysHire for a small number of tasks holds for a longer period of time.

In all plots of Figure~\ref{fig:results-problem-no-salary} we use coherence parameter $p=100$, which means we expect the input stream to be composed, on average, of runs of $100$ similar tasks (i.e., having Jaccard coefficient of at most 0.5 between consecutive ones). 
In this setting, even if the workload is not coherent (experimentally, even for $p=1$),  \algouseforfree is still better than \algoAlwaysOutsource.

\subsection{Experiments with \problemsalary}\label{subsec:results-tfo}

\spara{Baselines.}
As in \problemuseforfree, we consider baselines \algoAlwaysHire and \algoAlwaysOutsource. 
Additionally, we consider \naivetfo (defined in Section~\ref{concept:naivetfo}), which does not have a theoretical guarantee.
%

\spara{Results.}
Figure~\ref{fig:results-problem-salary} (Right) summarizes our results for \problemsalary.
We observe that \algowithsalary,\naivetfo, and \metaalgo have the smallest total cost,
followed by \algoAlwaysOutsource. In contrast, the \algoAlwaysHire strategy has
much higher cost due to mounting salary costs.
We also observe that while \naivetfo does not offer the theoretical guarantees of
\algowithsalary, it performs well in practice.

\begin{figure}[t]
\subfigure[UpWork: \algowithsalary vs. \algoAlwaysOutsource]{\includegraphics[width=.23\textwidth]{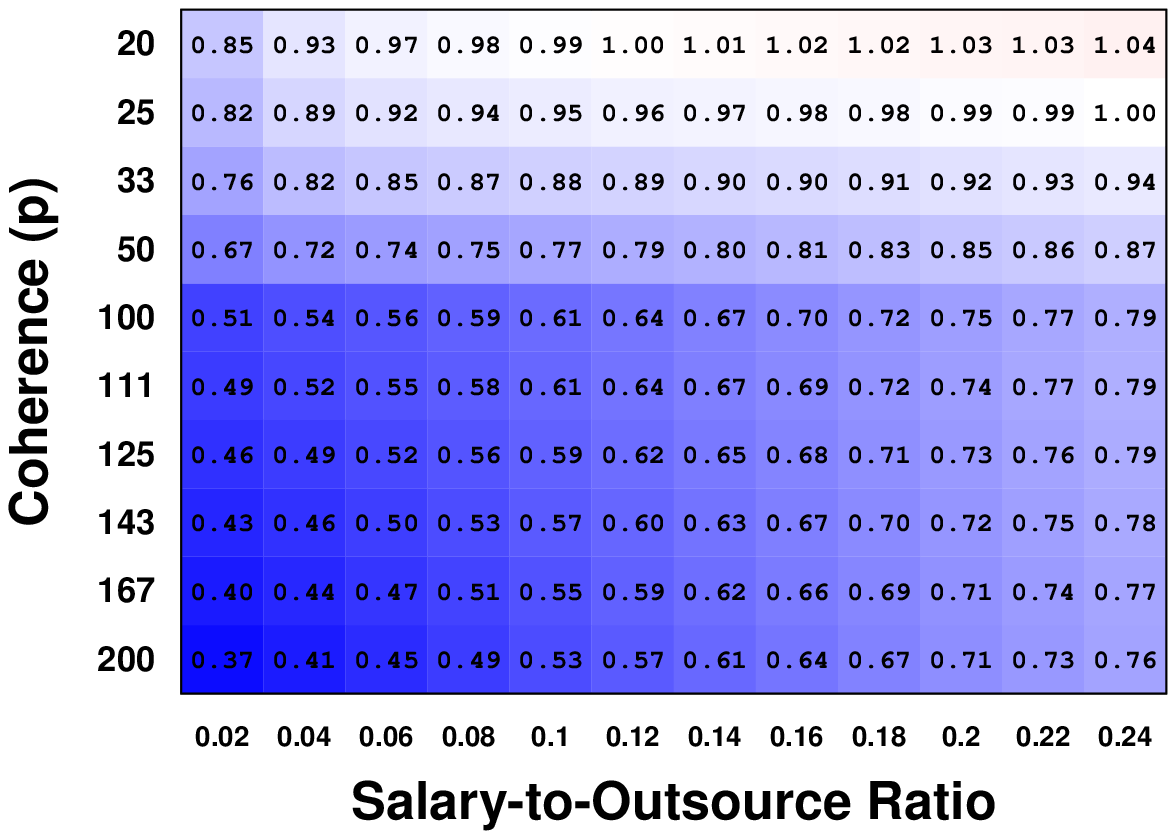}\label{fig:upwork-exploration}}
\subfigure[\mbox{UpWork: \metaalgo vs.} \algoAlwaysOutsource]{\includegraphics[width=.23\textwidth]{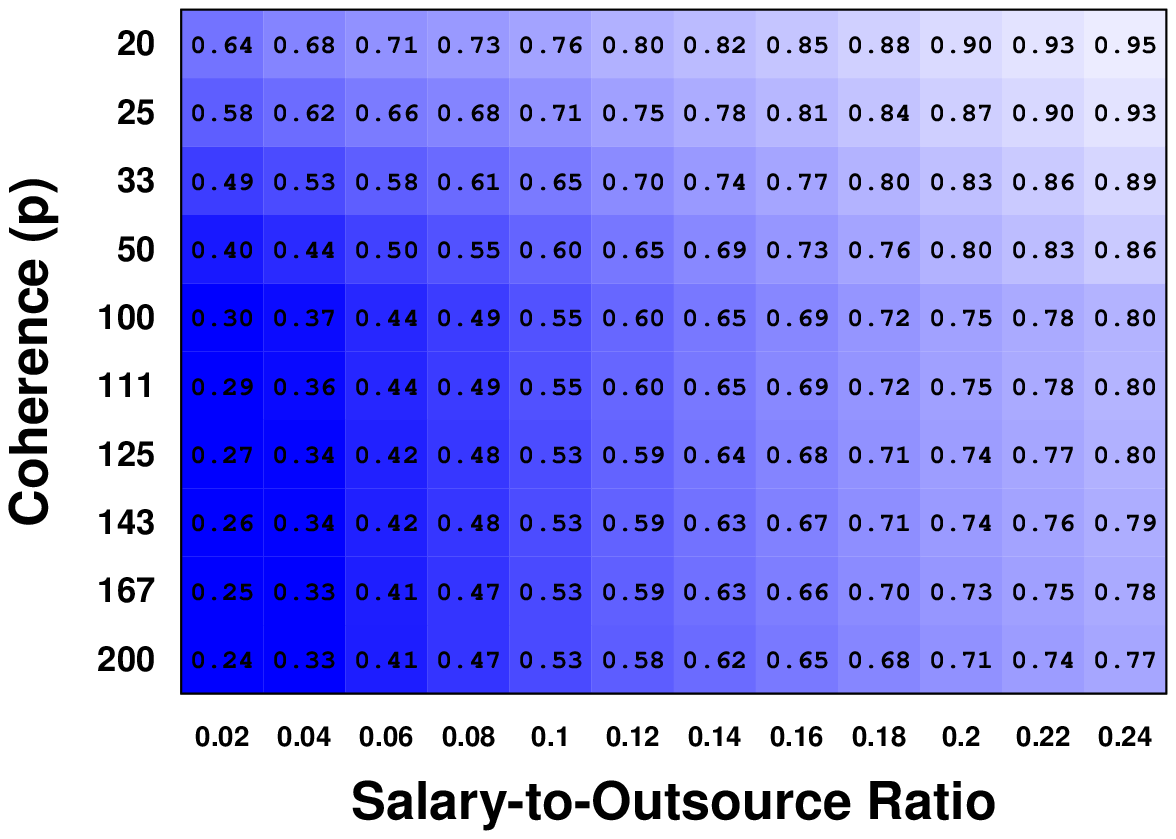}\label{fig:upwork-exploration-meta}}
\centering

\subfigure[\mbox{Freelancer: \algowithsalary vs.} \algoAlwaysOutsource]{\includegraphics[width=.23\textwidth]
{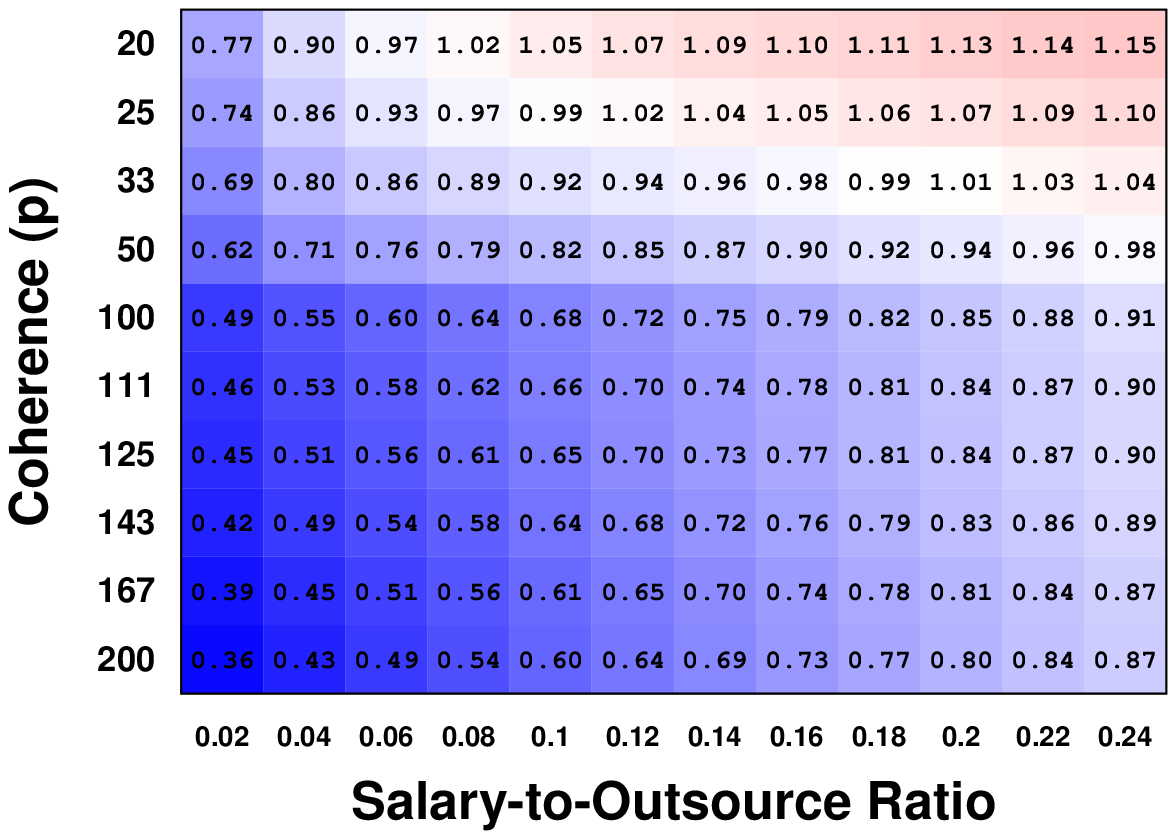}\label{fig:freelancer-exploration}}
\subfigure[\mbox{Freelancer: \metaalgo vs.} \algoAlwaysOutsource]{\includegraphics[width=.23\textwidth]
{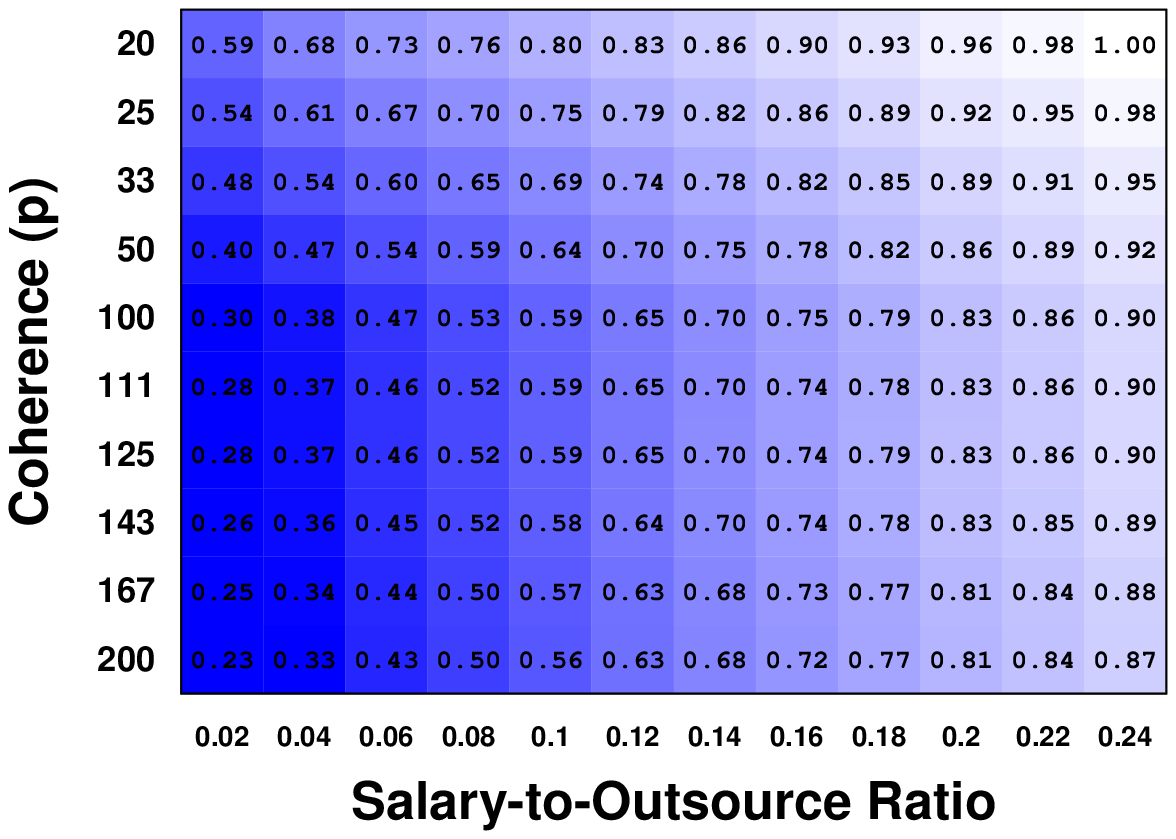}\label{fig:freelancer-exploration-meta}}
\centering

\subfigure[Guru: \algowithsalary vs. \algoAlwaysOutsource]{\includegraphics[width=.23\textwidth]{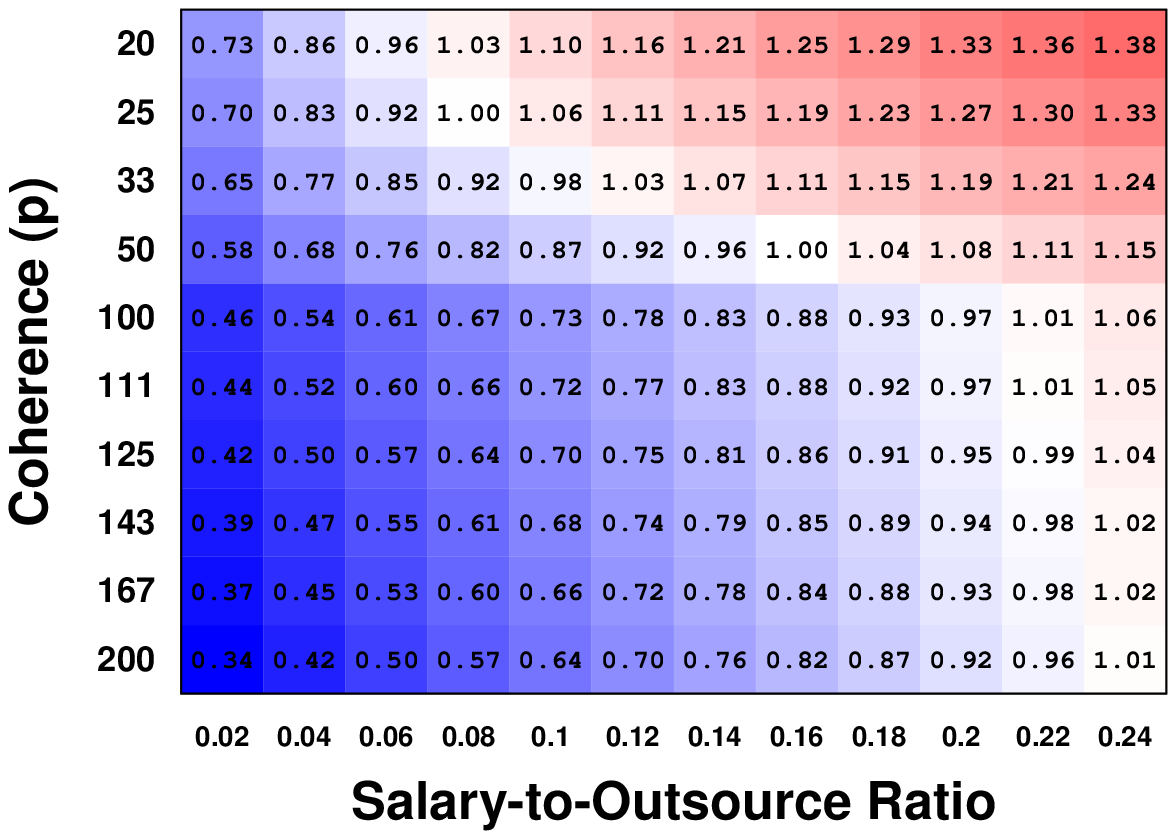}\label{fig:guru-exploration-exploration}}
\subfigure[\mbox{Guru: \metaalgo vs.} \algoAlwaysOutsource]{\includegraphics[width=.23\textwidth]{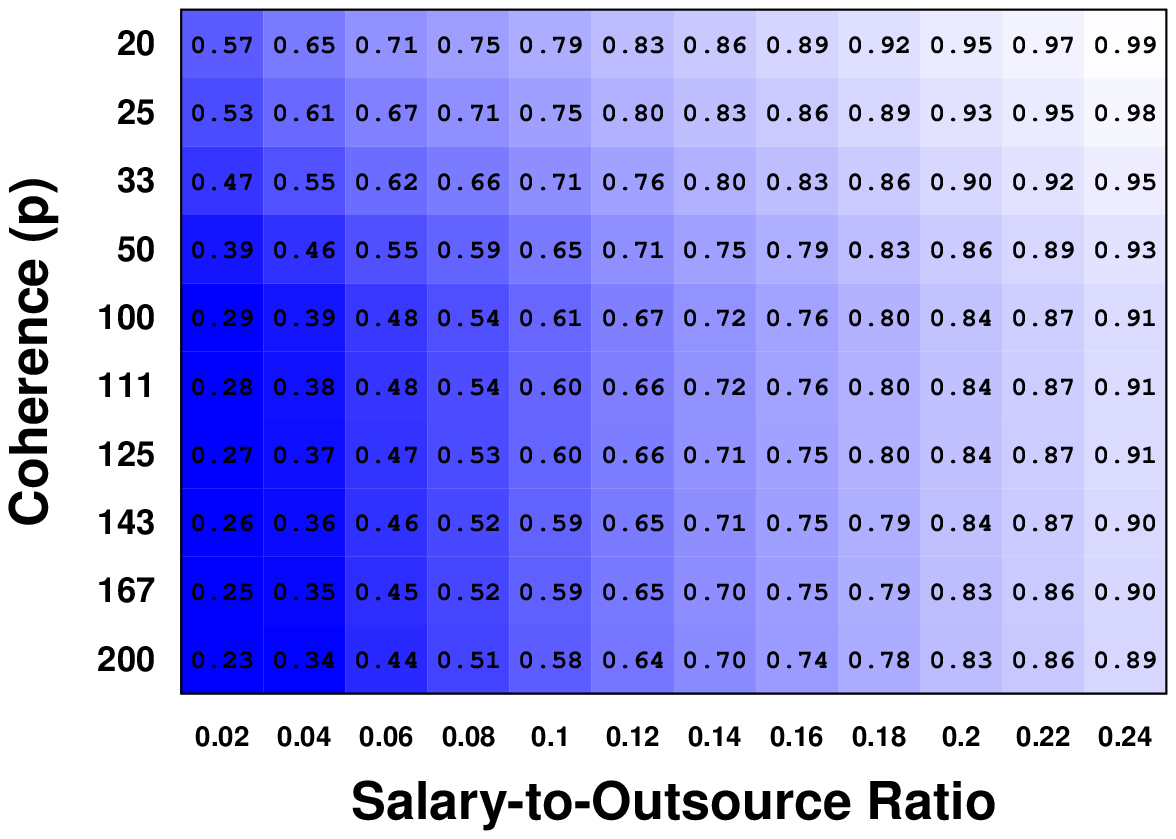}\label{fig:guru-exploration-meta}}
\caption{(Best seen in color.) Left: ratio of the cost achieved by \algowithsalary and \algoAlwaysOutsource. Right: ratio of the cost achieved by \metaalgo and \algoAlwaysOutsource. Coherence parameter $p$ varies from 0.02 to 0.24; salary-to-outsource ratio varies from 20 to 200; the number of tasks is 10K.
Colors represent the ratio of costs: blue (dominant towards the bottom-left) indicates the region where our algorithms \algowithsalary and \metaalgo have smaller cost, while red indicates the region where the baseline \algoAlwaysOutsource has smaller cost. In the white region both algorithms have similar costs.
}
\label{fig:exp-local-details}
\end{figure}

\spara{Variations.}
Similarly to \problemuseforfree, varying $\hiringcost{r}$ does not bring
dramatic changes, but as $\hiringcost{r}$ increases while maintaining
workload coherence and salary to outsource cost ratios constant, the
advantage of \algowithsalary over \algoAlwaysOutsource decreases, and
for large hiring costs \algoAlwaysOutsource has the smallest cost (plots
omitted for brevity).
Concretely, for $p=100$ and $\salary{r} = \outsourcecost{r}/10$,
if we vary the hiring cost $\hiringcost{r}$ (from $1\outsourcecost{r}$ to $30 \outsourcecost{r}$),
the total cost of \algowithsalary remains less or equal 
than the total cost of \algoAlwaysOutsource until $\hiringcost{r} = 16 \outsourcecost{r}$, when the cost of \algowithsalary becomes larger than the cost of \algoAlwaysOutsource for the workload generated using the Guru dataset.  The corresponding values of $\hiringcost{r}$ for workloads generated with Freelancer and Upwork data are $\hiringcost{r}=18\outsourcecost{r}$
and $\hiringcost{r}=26\outsourcecost{r}$ respectively.
As expected, if the hiring costs are sufficiently large, \algoAlwaysOutsource becomes a dominant strategy.
%

Figure~\ref{fig:exp-local-details} (Left) compares 
\algowithsalary and \algoAlwaysOutsource experimentally by varying the coherence
parameter $p$ from $20$ to $200$ and $\sigma_r$ from $\lambda_r / 50$ to
$\lambda_r / 4$. 
We observe that less coherent workloads and high salaries make hiring more expensive; \algoAlwaysOutsource then becomes a dominant strategy.
Figure~\ref{fig:exp-local-details} (Right) shows the power of the
\metaalgo algorithm. We observe that it performs equal or better than
\algoAlwaysOutsource for all the range of parameters.

\spara{Performance.}
Our code, which will be released with this paper, is a relatively straightforward mapping of the algorithm to simple counters. Written in Java, it requires about 5 to 8 seconds on average to process 10K incoming tasks using
commodity hardware.
We remark that, although our formulation is a linear program, the method does not involve \emph{solving} the linear program, 
instead, we obtain the solution using the specific primal--dual method that we have described and analyzed.


%% file: related.tex
To the best of our knowledge, 
we are the first to introduce and solve the Team Formation with Outsourcing (\ourproblem) problem.
However, our work is related to existing work on crowdsourcing, team formation, and online algorithms design, which we outline next.

\spara{Crowdsourcing.}
%
%
%
%
Among the extensive literature in crowdsourcing, the most related to ours is the work of~\citet{ho2012online}. Their goal is to assign individual workers to 
tasks, based on the workers' skills.
Although Ho and Vaughan also deploy the primal--dual 
technique
to solve the task-assignment problem,
the tasks they consider
can be performed by individual workers and not by teams.
Thus, both their problem and 
their algorithm is different from ours. 
 
\spara{Team formation.} 
A large body of work in team formation considers the following problem: given a social or a collaboration
network among the workers and a set of 
skills that needed to be covered, select a team of experts that can collectively cover all the required skills, while minimizing the communication cost
between the team members~\cite{anagnostopoulos12online, an13finding,dorn10composing,gajewar12multiskill,kargar11teamexp,lappas2009finding,li10team,sozio10community}.
Other variants of this problem have also considered
optimizing the cost of recruiting promising candidates for a set of pre-defined tasks in an offline fashion~\cite{golshan14profit} and minimizing the workload assigned to
each individual team member~\cite{majumder12capacitated,anagnostopoulos2010power}. 

Although the concept of set-cover is common between our work and previous work,
the framework we propose on this paper is different in multiple dimensions.
First,  we do not focus on optimizing the communication cost; 
in fact we do not assume any network among the individual workers. Our goal is to minimize
the overall cost paid on hiring, outsourcing, and salary costs. 
This difference in the objectives leads to different
(and new) optimization problems that we need to solve.
Secondly, most of the work above focuses on the offline version
of the team-formation problem, where the tasks to be completed are 
a-priori known to the algorithm. The exception is the work of \citet{anagnostopoulos2010power,anagnostopoulos12online}. 
%
 However, in their setting they aim to distribute the workload as evenly as possible among the workers,
while our objective is to minimize the overall cost of maintaining a team that can complete the 
arriving tasks. 
Moreover, the option of outsourcing that we propose is new with respect to the team formation literature. Finally, in the design
of our online algorithms we use the primal--dual framework, which was not the case for previous work on online team formation.

\spara{Primal--dual algorithms for online problems.}
The algorithms we design for our problems use the primal--dual technique. A thorough analysis on the applicability of this technique for online problems can be found in the book by
Buchbinder and Naor~\cite{buchbinder2009design} and in \cite{Bansal13}.
Probably the most closely related to problem are the \emph{ski-rental} and the
\emph{set cover}
problems. We have already discussed the connection of {\problemsalary} to ski-rental and set cover 
in Section~\ref{sec:preliminaries}.
One can also draw the analogy with caching;  one can think that bringing a page to the
main memory is analogous to hiring a person. The main differences are
that in the typical caching problem we do not have covering constraints, 
there are no recurring costs for keeping pages in the cache,
and there is a fixed limit on the number of pages we can insert in the
cache.

%% file: conclusions.tex
In this paper, we introduced and studied \emph{Team Formation with Outsourcing}. 
We showed that hiring, firing, and outsourcing decisions can be taken by an online algorithm 
%
leading to cost savings with respect to alternatives.
These cost savings are more striking when
\begin{inparaenum}[(1)]
 \item the hiring and salary costs are low, because then hiring becomes an attractive option;
 \item the tasks exhibit high coherence, i.e., consecutive tasks are similar to
each other; and 
 \item the time horizon is long enough that we can find a
core pool of workers to stay hired and satisfy a large fraction of the
skills required by incoming tasks.
\end{inparaenum}

Technically, the problems we have analyzed in this paper involve
embedding a set-cover problem in an online algorithm.
%
Our main algorithms (\algouseforfree, \algowithsalary) are able to give
results that are competitive in practice and, equally importantly,
theoretically close to the best one can hope for. The design of our
algorithms is based on the online primal--dual technique; we provide an
experimental evidence of the goodness of this method even for a complex 
real-world problem. 
Furthermore, we present two 
heuristics which,
although in theory are not competitive, perform
well in practice.
Future work may extend this by considering worker compatibility~\cite{lappas2009finding,anagnostopoulos12online},
learning of new skills by hired workers, or other extensions.

\spara{Future work.}
As most problems, we can introduce further elements to introduce even more generality.
For instance, the algorithms we have described assume one and only one
task arrives per unit of time, can be extended trivially to cases in
which task arrivals occur at arbitrary times.

As we noted in Section~\ref{sec:relwork}, there are also parallels with
scenarios of caching and paging. Extending \ourproblem when the number
of hired workers is limited turned out to be a challenging combination
of set cover, weighted caching and ski rental.
We have began to study these problems, our preliminary results show that
we can achieve a $O(\log k \log m)$ approximation, in which $k$ is the
maximum size of the worker pool.
A more natural constraint could be that,
for instance, the total cost paid per unit of time cannot exceed a
certain budget, which would represent a cap in weekly or monthly
personel expenses.
Another element we could incorporate is the possibility of not handling
a task,
but instead paying
a penalty when a task is too difficult to handle with current workers
and it is expensive to replace the worker pool with new workers. Other
variants can include workers with different ability levels.
We plan to study some of these variants in future work.

Additionally, we note that all the algorithms we have presented in this
paper are deterministic. Just as randomized algorithms for paging can be
defined in the primal--dual framework~\cite{buchbinder2009design}, it is
of interest to introduce other update rules for the primal variables that allow
us to describe a randomized algorithm.


\spara{Reproducibility.}
The code and data of this paper can be found at
\textsl{https://github.com/adrianfaz/Algorithms-for-Hiring-and-Outsourcing-in-the-Online-Labor-Market}.